\documentclass[12pt,letterpaper]{amsart}

\newcommand{\indentalign}{\hspace{0.3in}&\hspace{-0.3in}}
\newcommand{\la}{\langle}
\newcommand{\ra}{\rangle}
\newcommand{\R}{\text{Re}\,}

\newcommand{\ds}{\displaystyle}
\newcommand{\sgn}{\text{sgn}\,}
\newcommand{\supp}{\text{supp}\,}

\usepackage{amsmath}
\usepackage{amssymb}
\usepackage{amsthm}
\usepackage{amsxtra}

\newtheorem{theorem}{Theorem}
\newtheorem{definition}[theorem]{Definition}
\newtheorem{remark}[theorem]{Remark}

\newtheorem{lemma}[theorem]{Lemma}

\numberwithin{equation}{section}
\numberwithin{theorem}{section}

\title[IBVP for KdV]{The initial-boundary value problem for the Korteweg-de Vries equation}
\author{Justin Holmer}
\thanks{The content of this article appears as part of the author's Ph.D. thesis at the University of Chicago under the direction of Carlos Kenig.  The author is partially supported by an NSF postdoctoral fellowship.}
\address{University of California, Berkeley}


\begin{document}

\maketitle

\begin{abstract}
We prove local well-posedness of the initial-boundary value problem for the Korteweg-de Vries equation on right half-line, left half-line, and line segment, in the low regularity setting.  This is accomplished by introducing an analytic family of boundary forcing operators.
\end{abstract}

\maketitle

\tableofcontents

\newpage

\section{Introduction}
We shall study the following formulations of the initial-boundary value problem for the Korteweg-de Vries (KdV) equation.  On the right half-line $\mathbb{R}^+=(0,+\infty)$, we consider
\begin{equation} \label{SE:100}
\left\{
\begin{aligned}
&\partial_tu + \partial_x^3 u + u\partial_xu = 0  && \text{for }(x,t)\in (0,+\infty)\times (0,T) \\
& u(0,t) = f(t) && \text{for }t\in (0,T)\\
&u(x,0) = \phi(x) && \text{for }x\in (0,+\infty)
\end{aligned}
\right.
\end{equation}
On the left half-line $\mathbb{R}^-=(-\infty,0)$, we consider
\begin{equation} \label{SE:101}
\left\{
\begin{aligned}
&\partial_tu + \partial_x^3 u + u\partial_xu = 0  && \text{for }(x,t)\in (-\infty,0)\times (0,T) \\
& u(0,t) = g_1(t) && \text{for }t\in (0,T)\\
&\partial_xu(0,t) = g_2(t) && \text{for }t\in (0,T)\\
&u(x,0) = \phi(x) && \text{for }x\in (-\infty,0)
\end{aligned}
\right.
\end{equation}
The presence of one boundary condition in the right half-line problem \eqref{SE:100} versus two boundary conditions in the left half-line problem \eqref{SE:101} can be motivated by uniqueness calculations for smooth decaying solutions to the linear equation $\partial_t u + \partial_x^3u=0$.  Indeed, for such $u$ and $T>0$, we have
\begin{equation} \label{E:1}
\int_{x=0}^{+\infty} u(x,T)^2 \, dx= 
\begin{aligned}[t]
&\int_{x=0}^{+\infty} u(x,0)^2\, dx \\
&+2\int_{t=0}^T (u(0,t) \partial_x^2 u(0,t) - \partial_xu(0,t)^2) \, dt
\end{aligned}
\end{equation}
and
\begin{equation} \label{E:2}
\int_{x=-\infty}^0 u(x,T)^2 \, dx = 
\begin{aligned}[t]
&\int_{x=-\infty}^0 u(x,0)^2 \, dx\\
&-2\int_{t=0}^T (u(0,t) \partial_x^2 u(0,t) + \partial_xu(0,t)^2) \, dt .
\end{aligned}
\end{equation}
Assuming $u(x,0)=0$ for $x>0$ and $u(0,t)=0$ for $0<t<T$, we can conclude from \eqref{E:1} that $u(x,T)=0$ for $x>0$.  However, the existence of $u(x,t)\neq 0$ for $x<0$ such that $u(x,0)=0$ for $x<0$ and $u(0,t)=0$ for $0<t<T$ is not precluded by \eqref{E:2}.  In fact, such nonzero solutions do exist (see \S \ref{S:linear}). On the other hand, \eqref{E:2} does show that assuming $u(x,0)=0$ for $x<0$, $u(0,t)=0$ for $0<t<T$, and $\partial_xu(0,t)=0$ for $0<t<T$ forces $u(x,t)=0$ for $x<0$, $0<t<T$.  These uniqueness considerations carry over to the nonlinear equation $\partial_t u + \partial_x^3 u + u\partial_xu =0$, at least in the high regularity setting.  

Given the formulations \eqref{SE:100} and \eqref{SE:101}, it is natural to consider the following configuration for the line segment $0<x<L$ problem:
\begin{equation} \label{SE:102}
\left\{
\begin{aligned}
&\partial_tu + \partial_x^3 u + u\partial_xu = 0  && \text{for }(x,t)\in (0,L)\times (0,T) \\
& u(0,t) = f(t) && \text{for }t\in (0,T)\\
& u(L,t) = g_1(t) && \text{for }t\in (0,T)\\
&\partial_xu(L,t) = g_2(t) && \text{for }t\in (0,T)\\
&u(x,0) = \phi(x) && \text{for }x\in (0,L)
\end{aligned}
\right.
\end{equation}

Now we discuss appropriate spaces for the initial and boundary data, again examining the behavior of solutions to the linear problem on $\mathbb{R}$ for motivation.  On $\mathbb{R}$, we define the $L^2$-based inhomogeneous Sobolev spaces $H^s=H^s(\mathbb{R})$ by the norm $\| \phi \|_{H^s} = \| \la \xi \ra^s \hat{\phi}(\xi) \|_{L^2_\xi}$, where $\la \xi \ra = (1+|\xi|^2)^{1/2}$.  Let $e^{-t\partial_x^3}$ denote the linear homogeneous solution group on $\mathbb{R}$, defined by
\begin{equation} \label{E:3}
e^{-t\partial_x^3}\phi(x) = \tfrac{1}{2\pi}\int_\xi e^{it\xi^3} \hat{\phi}(\xi) \, d\xi ,
\end{equation}
so that $(\partial_t + \partial_x^3)e^{-t\partial_x^3}\phi(x)=0$ and $e^{-t\partial_x^3}\phi(x)\big|_{t=0} = \phi(x)$.
The \textit{local smoothing} inequalities of \cite{KPV91} for the operator \eqref{E:3} are
\begin{align*} 
\|\theta(t)e^{-t\partial_x^3}\phi \|_{L_x^\infty H_t^\frac{s+1}{3}} &\leq c \|\phi\|_{H^s}\\
\|\theta(t)\partial_x e^{-t\partial_x^3}\phi \|_{L_x^\infty H_t^\frac{s}{3}} &\leq c \|\phi\|_{H^s},
\end{align*}
which can be deduced directly from the definition \eqref{E:3} by a change of variable.  These are sharp in the sense that the Sobolev exponents $\frac{s+1}{3}$ and $\frac{s}{3}$ cannot be replaced by higher numbers.  In \S \ref{S:N}, we shall define analogues of the inhomogeneous Sobolev spaces on the half-line, $H^s(\mathbb{R}^+)$, $H^s(\mathbb{R}^-)$, and on the line segment, $H^s(0,L)$.  We are thus motivated to consider initial-boundary data pairs $(\phi, f) \in H^s(\mathbb{R}^+) \times H^\frac{s+1}{3}(\mathbb{R}^+)$ for \eqref{SE:100}, $(\phi, g_1, g_2) \in H^s(\mathbb{R}^-)\times  H^\frac{s+1}{3}(\mathbb{R}^+) \times H^\frac{s}{3}(\mathbb{R}^+)$ for \eqref{SE:101}, and $(\phi, f, g_1, g_2) \in H^s(0,L)\times  H^\frac{s+1}{3}(\mathbb{R}^+)\times  H^\frac{s+1}{3}(\mathbb{R}^+) \times H^\frac{s}{3}(\mathbb{R}^+)$ for \eqref{SE:102}.  From these motivations, we are inclined to consider this configuration optimal in the scale of $L^2$-based Sobolev spaces. 

Local well-posedness (LWP), i.e.\ existence, uniqueness, and uniform continuity of the data-to-solution map,  of the initial-value problem (IVP) 
\begin{equation} \label{E:4}
\left\{
\begin{aligned}
&\partial_tu + \partial_x^3 u + u\partial_x u = 0 && \text{for }(x,t)\in \mathbb{R}\times \mathbb{R}\\
& u(x,0) = \phi(x) && \text{for }(x,t)\in \mathbb{R}
\end{aligned}
\right.
\end{equation}
has been studied by a number of authors over the past three decades.   For $s>\frac{3}{2}$, an \textit{a priori} bound can be obtained by the energy method and a solution can be constructed via the artificial viscosity method.  To progress to rougher spaces, it is necessary to invoke techniques of harmonic analysis to quantitatively capture the dispersion of higher frequency waves.  For $s>\frac{3}{4}$, \cite{KPV91} proved LWP of \eqref{E:4} by the contraction method in a space built out of various space-time norms, using oscillatory integral and local smoothing estimates. For $s>-\frac{3}{4}$, \cite{B93} \cite{KPV94} \cite{KPV96} proved LWP of \eqref{E:4} via the contraction method in Bourgain spaces (denoted in the literature as $X_{s,b}$), which are constructed to delicately analyze the interaction of waves in different frequency zones.  LWP for $s=-\frac{3}{4}$ is proved in \cite{CCT03} by using the Miura transform to convert KdV to mKdV (nonlinearity $u^2\partial_xu$) where the corresponding endpoint result is known.  These authors also prove local ill-posedness of \eqref{E:4} for $s<-\frac{3}{4}$ in the sense that the data-to-solution map fails to be uniformly continuous.  If one only requires that the data-to-solution map be continuous ($C^0$ well-posedness), and not uniformly continuous, then the regularity requirements can possibly be relaxed further.  Although this has not yet been shown for the KdV equation on the line, \cite{KT03} have proved, for the KdV equation on the circle $\mathbb{T}$,  $C^0$ local well-posedness in $H^{-1}(\mathbb{T})$, whereas it has been shown by \cite{CCT03} that the data-to-solution map cannot be uniformly continuous in $H^s(\mathbb{T})$ for $s<-\frac{1}{2}$.

Our goal in studying \eqref{SE:100} is to obtain low regularity results.  It therefore seems reasonable to restrict to $-\frac{3}{4}<s<\frac{3}{2}$.  We shall omit $s=\frac{1}{2}$ due to difficulties in formulating the compatibility condition (see below).  A Dini integral type compatibility condition would probably suffice at this point, although we have decided not to explore it.  We have also decided not to explore the case $s=-\frac{3}{4}$ or the likely ill-posedness result for \eqref{SE:100} and \eqref{SE:101} when $s<-\frac{3}{4}$.  

Note that the trace map $\phi \to \phi(0)$ is well-defined on $H^s(\mathbb{R}^+)$ when $s>\frac{1}{2}$.  If $s>\frac{1}{2}$, then $\frac{s+1}{3}>\frac{1}{2}$, and both $\phi(0)$ and $f(0)$ are well-defined quantities.  Since $\phi(0)$ and $f(0)$ are both meant to represent $u(0,0)$, they must agree.  On the other hand, if $s<\frac{3}{2}$, then $s-1<\frac{1}{2}$ and $\frac{s}{3}<\frac{1}{2}$, so in \eqref{SE:101}, neither $\partial_x u \in H^{s-1}$ nor $g_2\in H^\frac{s}{3}$ have a well-defined trace at $0$.  

Therefore, we consider \eqref{SE:100} for $-\frac{3}{4}< s<\frac{3}{2}$, $s\neq \frac{1}{2}$ in the setting
\begin{equation} \label{E:111}
\phi\in H^s(\mathbb{R}^+), \; f\in H^\frac{s+1}{3}(\mathbb{R}^+), \; \text{and if }\tfrac{1}{2}<s<\tfrac{3}{2}, \; \phi(0)=f(0) .
\end{equation}
We consider \eqref{SE:101} for $-\frac{3}{4}< s<\frac{3}{2}$, $s\neq \frac{1}{2}$ in the setting
\begin{equation} \label{E:112}
\begin{gathered} 
\phi\in H^s(\mathbb{R}^-), \; g_1\in H^\frac{s+1}{3}(\mathbb{R}^+),  \; g_2\in H^\frac{s}{3}(\mathbb{R}^+) \\
\text{and if }\tfrac{1}{2}<s<\tfrac{3}{2}, \; \phi(0)=g_1(0)
\end{gathered} 
\end{equation}
We consider \eqref{SE:102} for $-\frac{3}{4}< s<\frac{3}{2}$, $s\neq \frac{1}{2}$ in the setting
\begin{equation}\label{E:113}
\begin{gathered} 
\phi\in H^s(0,L), \; f\in H^\frac{s+1}{3}(\mathbb{R}^+),  \; g_1\in H^\frac{s+1}{3}(\mathbb{R}^+),  \; g_2\in H^\frac{s}{3}(\mathbb{R}^+) \\
 \text{and if }\tfrac{1}{2}<s<\tfrac{3}{2}, \; \phi(0)=f(0), \phi(L)=g_1(0)
\end{gathered}
\end{equation}

The solutions we construct shall have the following characteristics.
\begin{definition} \label{D:strong}
$u(x,t)$ will be called a \emph{distributional solution of \eqref{SE:100}, \eqref{E:111} [resp. \eqref{SE:101}, \eqref{E:112}] on $[0,T]$} if
\begin{enumerate}
\item \emph{Well-defined nonlinearity}:  $u$ belongs to some space $X$ with the property that $u\in X \Longrightarrow \partial_x u^2$ is a well-defined distribution.
\item $u(x,t)$ satisfies the equation \eqref{SE:100} [resp. \eqref{SE:101}] in the sense of distributions on the set $(x,t) \in (0,+\infty) \times (0,T)$ [resp. $(x,t) \in (-\infty,0) \times (0,T)$].
\item \emph{Space traces:}  $u\in C( [0,T]; \; H^s_x)$ and in this sense $u(\cdot,0) = \phi$ in  $H^s(\mathbb{R}^+)$ [resp. $u(\cdot,0) = \phi$ in  $H^s(\mathbb{R}^-)$].
\item \emph{Time traces:} $u\in C( \mathbb{R}_x ; H^\frac{s+1}{3}(0,T))$ and in this sense $u(0, \cdot )=f$ in $H^\frac{s+1}{3}(0,T)$ [resp. $u(0, \cdot )=g_1$ in $H^\frac{s+1}{3}(0,T)$].
\item \emph{Derivative time traces:} $\partial_x u\in C( \mathbb{R}_x ; H^\frac{s}{3}(0,T))$ and only for \eqref{SE:101},\eqref{E:112} we require that in this sense, $u(0, \cdot )=g_2$ in $H^\frac{s}{3}(0,T)$.
\end{enumerate}
\end{definition}
In our case, $X$ shall be the modified Bourgain space $X_{s,b}\cap D_\alpha$ with $b<\frac{1}{2}$ and $\alpha>\frac{1}{2}$, where
\begin{equation}
\label{E:400}
\begin{gathered}
\|u\|_{X_{s,b}} = \left( \iint_{\xi,\tau} \la\xi\ra^{2s}\la\tau-\xi^3\ra^{2b}|\hat{u}(\xi,\tau)|^2 \, d\xi \, d\tau \right)^{1/2},\\
\|u\|_{D_\alpha} =  \left( \iint_{|\xi|\leq 1} \la\tau\ra^{2\alpha} |\hat{u}(\xi,\tau)|^2 \, d\xi \,d \tau \right)^{1/2}.
\end{gathered}
\end{equation}
The space $X_{s,b}$, with $b>\frac{1}{2}$, is typically employed in the study of the IVP \eqref{E:4}.  For $b>\frac{1}{2}$, the bilinear estimate (Lemma \ref{L:bilinear}) holds without the low frequency modification $D_\alpha$, and thus $D_\alpha$ is not necessary in the study of the IVP.  The introduction of the Duhamel boundary forcing operator in our study of the IBVP, however, forces us to take $b<\frac{1}{2}$, and then $D_\alpha$ must be added in order for Lemma \ref{L:bilinear} to hold.  

A definition for \eqref{SE:102}, \eqref{E:113} can be given in the obvious manner.  We shall next introduce the concept of mild solution used by \cite{BSZ04}.

\begin{definition}
$u(x,t)$ is a \emph{mild solution} of \eqref{SE:100} [resp. \eqref{SE:101}] on $[0,T]$ if $\exists$ a sequence $\{  u_n  \}$ in $C( [0,T]; \; H^3(\mathbb{R}_x^+) )\cap C^1([0,T]; \; L^2(\mathbb{R}_x^+) )$ such that 
\begin{enumerate}
\item $u_n(x,t)$ solves \eqref{SE:100} in $L^2(\mathbb{R}_x^+)$ [resp. \eqref{SE:101} in $L^2(\mathbb{R}_x^-)$] for $0<t<T$.
\item $\displaystyle \lim_{n\to +\infty} \|u_n -u \|_{C([0,T];\, H^s(\mathbb{R}_x^+))} =0$ [resp. $\displaystyle \lim_{n\to +\infty} \|u_n -u \|_{C([0,T];\, H^s(\mathbb{R}_x^-))} =0$].
\item $\displaystyle \lim_{n\to +\infty} \|u_n(0,\cdot)-f\|_{H^\frac{s+1}{3}(0,T)} =0$ [resp. $\displaystyle \lim_{n\to +\infty} \|u_n(0,\cdot)-g_1\|_{H^\frac{s+1}{3}(0,T)} =0$, $\displaystyle \lim_{n\to +\infty} \|\partial_x u_n(0,\cdot)-g_2\|_{H^\frac{s}{3}(0,T)} =0$].
\end{enumerate}
\end{definition}

\cite{BSZ05} have recently introduced a method for proving uniqueness of mild solutions for \eqref{SE:100}, \eqref{E:111}.

Our main result is the following existence statement.
\begin{theorem} \label{T:main} Let $-\frac{3}{4}<s<\frac{3}{2}$, $s\neq \frac{1}{2}$.
\begin{enumerate} 
\item \label{I:right}Given $(\phi,f)$ satisfying \eqref{E:111}, $\exists \; T>0$ depending only on the norms of $\phi$, $f$ in \eqref{E:111} and $\exists\; u(x,t)$ that is both a mild and distributional solution to \eqref{SE:100}, \eqref{E:111} on $[0,T]$.
\item \label{I:left}Given $(\phi,g_1,g_2)$ satisfying \eqref{E:112}, $\exists \; T>0$ depending only on the norms of $\phi$, $g_1$, $g_2$ in \eqref{E:112} and $\exists\; u(x,t)$ that is both a mild and distributional solution to \eqref{SE:101}, \eqref{E:112} on $[0,T]$.
\item \label{I:lineseg}Given $(\phi,f, g_1,g_2)$ satisfying \eqref{E:113}, $\exists \; T>0$ depending only on the norms of $\phi$, $f$, $g_1$, $g_2$ in \eqref{E:113} and $\exists\; u(x,t)$ that is both a mild and distributional solution to \eqref{SE:102}, \eqref{E:113} on $[0,T]$.
\end{enumerate}
In each of the above cases, the data-to-solution map is analytic as a map from the spaces in \eqref{E:111}, \eqref{E:112}, \eqref{E:113} to the spaces in Definition \ref{D:strong}.
\end{theorem}

The proof of Theorem \ref{T:main} involves the introduction of an analytic family of boundary forcing operators extending the single operator introduced by \cite{CK02} (further comments in \S \ref{S:overview}).

The main new feature of our work is the low regularity requirements for $\phi$ and $f$.  Surveys of the literature are  given in  \cite{BSZ02} \cite{BSZ03} and \cite{CK02}.  Here, we briefly mention some of the more recent contributions.  The  problem \eqref{SE:102}\eqref{E:113} for $s\geq 0$ is treated in  \cite{BSZ03} and \eqref{SE:100} \eqref{E:111} for $s>\frac{3}{4}$ in \cite{BSZ02} by a Laplace transform technique.  In a preprint appearing after this paper was submitted, \cite{BSZ06} have shown LWP of the problem \eqref{SE:100} for $s>-1$ with $H^s(\mathbb{R}^+)$ in \eqref{E:111} replaced by the weighted space
$$H^s(\mathbb{R}^+) = \{ \, \phi \in H^s(\mathbb{R}^+) \; \mid \: e^{\nu x}\phi(x) \in H^s(\mathbb{R}^+) \, \}$$
for $\nu>0$.  They further show LWP of the problem \eqref{SE:102},\eqref{E:113} for $s>-1$, thus improving Theorem \ref{T:main}\eqref{I:lineseg}.  In both of these results, the data-to-solution map is analytic, in contrast to the results of \cite{KT03} mentioned above.  A global well-posedness result for the problem \eqref{SE:100}\eqref{E:111} is obtained by \cite{Fam03} for $s\geq 0$.  Inverse scattering techniques have been applied to the problem \eqref{SE:101} by \cite{Fok02} and the linear analogue of the problem \eqref{SE:102} in \cite{FP01b} for Schwartz class data. 

I have carried out similar results for the nonlinear Schr\"odinger equation \cite{Hol05}.

\noindent\textbf{Acknowledgements}.  I would like to thank my Ph.D.\ advisor Carlos Kenig for invaluable guidance on this project.  I would also like to thank the referee for a careful reading and helpful suggestions.

\section{Overview}
\label{S:overview}

In this section, after giving some needed preliminaries, we introduce the Duhamel boundary forcing operator of \cite{CK02} and first apply it and a related operator to solve linear versions of the problems \eqref{SE:100}, \eqref{SE:101}.  Then we explain the need for considering a more general class of operators to address the nonlinear versions in $H^s$ for $-\frac{3}{4}<s<\frac{3}{2}$, $s\neq \frac{1}{2}$.

Since precise numerical coefficients become important, let us set down the convention 
$$\hat{f}(\xi) = \int_x e^{-ix\xi} f(x)\, dx .$$
Also, define $C_0^\infty(\mathbb{R}^+)$ as those smooth functions on $\mathbb{R}$ with support contained in $[0,+\infty)$.  Let $C_{0,c}^\infty(\mathbb{R}^+) = C_0^\infty(\mathbb{R}^+)\cap C_c^\infty(\mathbb{R})$.
The tempered distribution $\frac{t_+^{\alpha-1}}{\Gamma(\alpha)}$ is defined as a locally integrable function for $\text{Re }\alpha>0$, i.e.\
$$\left< \frac{t_+^{\alpha-1}}{\Gamma(\alpha)}, f \right> = \frac{1}{\Gamma(\alpha)}\int_0^{+\infty} t^{\alpha-1} f(t) \, dt .$$
Integration by parts gives, for \text{Re }$\alpha>0$, that 
\begin{equation}
\label{E:rev1}
\frac{t_+^{\alpha-1}}{\Gamma(\alpha)} = \partial_t^k \left[ \frac{t_+^{\alpha+k-1}}{\Gamma(\alpha+k)} \right]
\end{equation}
for all $k\in \mathbb{N}$.
  This formula can be used to extend the definition (in the sense of distributions) of $\frac{t_+^{\alpha-1}}{\Gamma(\alpha)}$ to all $\alpha \in \mathbb{C}$.  In particular, we obtain
$$\left. \frac{t_+^{\alpha-1}}{\Gamma(\alpha)}\right|_{\alpha=0} = \delta_0(t) .$$
A change of contour calculation shows that 
\begin{equation}
\label{E:410}
\left[\frac{t_+^{\alpha-1}}{\Gamma(\alpha)} \right]\sphat(\tau)=e^{-\frac{1}{2}\pi i \alpha}(\tau - i0)^{-\alpha}
\end{equation}
where $(\tau-i0)^{-\alpha}$ is the distributional limit. If $f\in C_0^\infty(\mathbb{R}^+)$, we define
$$\mathcal{I}_\alpha f  = \frac{t_+^{\alpha-1}}{\Gamma(\alpha)} * f .$$
Thus, when $\text{Re }\alpha>0$,
$$\mathcal{I}_\alpha f(t)  = \frac{1}{\Gamma(\alpha)}\int_0^t (t-s)^{\alpha-1} f(s) \, ds$$
and $\mathcal{I}_0f=f$, $\mathcal{I}_1f(t)=\int_0^t f(s) \, ds$, and $\mathcal{I}_{-1}f=f'$.  Also $\mathcal{I}_{\alpha}\mathcal{I}_\beta = \mathcal{I}_{\alpha+\beta}$, which follows from \eqref{E:410}.  For further details on the distribution $\frac{t_+^{\alpha-1}}{\Gamma(\alpha)}$, see \cite{F98}. 
\begin{lemma} \label{L:RL}
If $f\in C_0^{\infty}(\mathbb{R}^+)$, then $\mathcal{I}_{\alpha}f \in C_0^{\infty}(\mathbb{R}^+)$, for all $\alpha\in \mathbb{C}$.
\end{lemma}
\begin{proof}
By \eqref{E:rev1} and integration by parts, it suffices to consider the case $\R \alpha>1$.  In this case, it is clear that $\supp I_\alpha f \subset [0,+\infty)$ and it remains only to show that $I_\alpha f(t)$ is smooth.  By a change of variable
$$I_\alpha f(t) = \frac{1}{\Gamma(\alpha)} \int_0^t s^{\alpha-1} f(t-s) \, ds .$$
Smoothness of $I_\alpha f(t)$ follows by the fundamental theorem of calculus, differentiation under the integral sign, and that $\partial_t^kf(0)=0$ for all $k$.
\end{proof}

The \textit{Airy function} is 
$$A(x) = \frac{1}{2\pi} \int_\xi e^{ix\xi}e^{i\xi^3} \, d\xi .$$
$A(x)$ is a smooth function with the asymptotic properties
\begin{align*}
A(x) &\sim c_1 x^{-1/4} e^{-c_2 x^{3/2}}(1+O(x^{-3/4})) && \text{as }x\to +\infty \\
A(-x) &\sim c_2 x^{-1/4}\cos(c_2x^{3/2}-\tfrac{\pi}{4})(1+O(x^{-3/4})) && \text{as }x\to +\infty
\end{align*}
for specific $c_1,c_2>0$ (see, e.g.\ \cite{SS04}, p.\ 328).  We shall below need the values of $A(0)$, $A'(0)$, and $\int_0^{+\infty} A(y)\, dy$, and so we now compute them. 
$$A(0)= \frac{1}{2\pi} \int_\xi e^{i\xi^3} \, d\xi = \frac{1}{6\pi} \int_\eta \eta^{-2/3}e^{i\eta} \, d\eta = \frac{\frac{\sqrt 3}{2}\Gamma(\frac{1}{3})}{3\pi} = \frac{1}{3\Gamma(\frac{2}{3})}$$
by a change of contour calculation, and in the final step, an application of the identity $\Gamma(z)\Gamma(1-z)=\pi/\sin \pi z$.  Similarly one finds
$$A'(0) = \frac{1}{2\pi} \int_\xi i\xi e^{i\xi^3} \, d\xi = -\frac{1}{3\Gamma(\frac{1}{3})} .$$
Also,
$$\int_{y=0}^{+\infty} A(y)\, dy = \frac{1}{2\pi} \int_\xi \int_{y=0}^{+\infty} e^{iy\xi} \, dy \; e^{i\xi^3} \, d\xi = \frac{1}{2\pi} \int_\xi \hat{H}(-\xi) e^{i\xi^3} \, d\xi$$
where $H(y)=0$ for $y<0$, $H(y)=1$ for $y>0$ is the Heaviside function. Now (see \cite{F98}, p.\ 101) $\hat{H}(\xi) = \text{p.v.}\frac{1}{i\xi} + \pi \delta_0(\xi)$, which inserted above and combined with the identity $(\text{p.v.}1/x)\sphat(\xi) = -i\pi\sgn \xi$ yields 
$$\int_0^{+\infty} A(y)\, dy = \frac{1}{3} .$$

\subsection{Linear versions}
\label{S:linear}

We define the Airy \textit{group} as 
\begin{equation}
\label{E:G}
e^{-t\partial_x^3}\phi(x) = \frac{1}{2\pi} \int_\xi e^{ix\xi} e^{it\xi^3} \hat{\phi}(\xi) \, d\xi
\end{equation}
so that 
\begin{equation}
\label{E:Geq}
\left\{
\begin{aligned}
& (\partial_t  + \partial_x^3) [e^{-t\partial_x^3}\phi](x,t)  = 0 && \text{for }(x,t)\in \mathbb{R}\times \mathbb{R}\\
& [e^{-t\partial_x^3}\phi](x,0) = \phi(x) && \text{for }x\in \mathbb{R}
\end{aligned}
\right.
\end{equation}

We now introduce the \textit{Duhamel boundary forcing operator} of \cite{CK02}.  For $f\in C_0^\infty(\mathbb{R}^+)$, let
\begin{equation}
\label{E:Dbf}
\begin{aligned}
\mathcal{L}^0f(x,t) &= 3 \int_0^t e^{-(t-t')\partial_x^3}\delta_0(x)\mathcal{I}_{-2/3}f(t')\, dt' \\
&= 3\int_0^t A\left( \frac{x}{(t-t')^{1/3}} \right) \frac{\mathcal{I}_{-2/3}f(t')}{(t-t')^{1/3}} \, dt'
\end{aligned}
\end{equation}
so that
\begin{equation}
\label{E:Dbfeq1}
\left\{
\begin{aligned}
& (\partial_t + \partial_x^3) \mathcal{L}^0f(x,t) = 3\delta_0(x)\mathcal{I}_{-2/3}f(t) && \text{for }(x,t)\in \mathbb{R}\times \mathbb{R}\\
& \mathcal{L}^0f(x,0) = 0 && \text{for }x\in \mathbb{R}\\
\end{aligned}
\right.
\end{equation}
We begin with the spatial continuity and decay properties of $\mathcal{L}^0f$, $\partial_x \mathcal{L}^0f$, and $\partial_x^2\mathcal{L}^0f$, for $f\in C_0^\infty(\mathbb{R}^+)$.

\begin{lemma} \label{L:Dbf2}
Let $f\in C_0^\infty(\mathbb{R}^+)$.  Then for fixed $0\leq t\leq 1$, $\mathcal{L}^0f(x,t)$ and $\partial_x\mathcal{L}^0f(x,t)$ are continuous in $x$ for all $x\in \mathbb{R}$ and satisfy the spatial decay bounds 
\begin{equation} \label{AE:463}
|\mathcal{L}^0f(x,t)| + |\partial_x \mathcal{L}^0f(x,t)| \leq c_k\|f\|_{H^{k+1}} \langle x \rangle^{-k} \quad \forall \; k\geq 0 .
\end{equation}
For fixed $0\leq t\leq 1$, $\partial_x^2\mathcal{L}^0f(x,t)$ is continuous in $x$ for $x\neq 0$ and has a step discontinuity of size $3\mathcal{I}_{2/3}f(t)$ at $x=0$. Also, $\partial_x^2\mathcal{L}^0f(x,t)$ satisfies the spatial decay bounds 
\begin{equation} \label{AE:2021}
|\partial_x^2 \mathcal{L}^0f(x,t)| \leq c_k \|f\|_{H^{k+2}} \la x \ra^{-k} \qquad  \forall \; k\geq 0
\end{equation}
\end{lemma}
\begin{proof}
To establish \eqref{AE:463}, it suffices to show that $\| \langle \xi \rangle \partial_\xi^k\,  \widehat{\mathcal{L}^0f}(\xi,t)\|_{L^1_\xi} \leq c_k\|f\|_{H^k}$, $\forall \; k\geq 0$.  Let $\phi(\xi,t) = \int_0^t e^{i(t-t')\xi}h(t') \, dt'$ for some (yet to be prescribed) $h\in C_0^\infty(\mathbb{R}^+)$.
We have
\begin{equation}\label{AE:2000}
\partial_\xi^k \phi(\xi,t) = i^k \int_0^t (t-t')^k  e^{i(t-t')\xi} h(t') \, dt' .
\end{equation}
By integration by parts in $t'$, 
\begin{equation} \label{AE:2001}
\partial_\xi^k\phi(\xi,t) = 
\begin{aligned}[t]
&\frac{i(-1)^{k+1}k!}{\xi^{k+1}} \int_0^t e^{i(t-t')\xi} \partial_{t'}h(t') \, dt' + \frac{i(-1)^kk!}{\xi^{k+1}}h(t)\\
& + \frac{i(-1)^{k+1}}{\xi^{k+1}} \int_0^t e^{i(t-t')\xi} \partial_{t'} \sum_{\substack{\alpha+\beta=k \\ \alpha\leq k-1}} c_{\alpha, \beta} \partial_{t'}^\alpha(t-t')^k \partial_{t'}^\beta h(t') \, dt'
\end{aligned}
\end{equation}
By \eqref{AE:2000}, \eqref{AE:2001} and the time localization, $|\partial_\xi^k\phi(\xi,t)| \leq c_k \|h\|_{H^k} \la \xi \ra^{-k-1}$.  Since $\widehat{\mathcal{L}^0f}(\xi,t) = \phi(\xi^3,t)$ with $h=3\mathcal{I}_{-2/3}f$, we have by Lemma \ref{L:RL1} that $|\partial_\xi^k\widehat{\mathcal{L}^0f}(\xi,t)|\leq c_k\|f\|_{H^{k+1}}\la \xi \ra^{-k-3}$, establishing \eqref{AE:463}.  By integration by parts in $t'$ in \eqref{E:Dbf},
\begin{equation} \label{AE:2006}
\partial_x^3 \mathcal{L}^0f(x,t) = 3\delta_0(x)\mathcal{I}_{-2/3}f(t) - \mathcal{L}^0(\partial_tf)(x,t) .
\end{equation}
This, together with the continuity properties of $\mathcal{L}^0(\partial_t f)$, shows that $\partial_x^2 \mathcal{L}^0f(x,t)$ is continuous in $x$ for $x\neq 0$ and has a step discontinuity of size $3\mathcal{I}_{-2/3}f(t)$ at $x=0$.  To see that $\partial_x^2\mathcal{L}^0f(x,t) \to 0$ as $x\to \pm \infty$, we first note that for $x<-1$, $\partial_x^2\mathcal{L}^0f(x,t) = \partial_x^2\mathcal{L}^0f(-1,t) - \int_x^{-1} \partial_y^3 \mathcal{L}^0f(y,t) \, dy$.  By \eqref{AE:2006} and \eqref{AE:463}, we can send $x\to -\infty$ and obtain that $\partial_x^2\mathcal{L}^0f(x,t) \to c$, for some constant $c$, as $x\to -\infty$.  Since $\partial_x \mathcal{L}^0f(0,t) = \int_{-\infty}^0 \partial_x^2\mathcal{L}^0f(y,t) \, dy$, we must have $c=0$.  We can similarly show that $\partial_x^2 \mathcal{L}^0f(x,t) \to 0$ as $x\to +\infty$.  For $x<0$, use $\partial_x^2 \mathcal{L}^0f(x,t) = \int_{-\infty}^x \partial_y^3 \mathcal{L}^0f(y,t) \, dy$, and for $x>0$, use $\partial_x^2\mathcal{L}^0f(x,t) = -\int_x^{+\infty} \partial_y^3\mathcal{L}^0f(y,t) \, dy $, together with \eqref{AE:463} and \eqref{AE:2006} to obtain the bound \eqref{AE:2021}.
\end{proof}

By Lemma \ref{L:Dbf2},  if $f\in C_0^\infty(\mathbb{R}^+)$, then $\mathcal{L}^0f(x,t)$ is continuous in $x$ on $\mathbb{R}$.  Since $A(0)=(3\Gamma(\frac{2}{3}))^{-1}$, the second representation of $\mathcal{L}^0f(x,t)$ in \eqref{E:Dbf} gives 
\begin{equation}
\label{E:Dbfeq2}
\mathcal{L}^0f(0,t) = f(t).
\end{equation}
It is thus clear that if we set
$$u(x,t) = e^{-t\partial_x^3}\phi(x) + \mathcal{L}^0\left( f- e^{-\cdot\partial_x^3}\phi\big|_{x=0} \right)(t)$$
then $u(x,t)$ solves the linear problem
$$\left\{
\begin{aligned}
&(\partial_t + \partial_x^3)u(x,t) = 0 && \text{for }x\neq 0 \\
&u(x,0) = \phi(x) &&\text{for }x\in \mathbb{R}\\
&u(0,t) = f(t) &&\text{for }t\in \mathbb{R}
\end{aligned}
\right.
$$
This would suffice, then, to solve the linear analogue of the right half-line problem \eqref{SE:100}, which has only one boundary condition.

Now we consider the linear analogue of the left half-line problem \eqref{SE:101}, which has two boundary conditions. 
Consider, in addition to $\mathcal{L}^0$, the second boundary forcing operator
\begin{equation}
\label{E:Dbf2}
\begin{aligned}
\mathcal{L}^{-1}f(x,t) &= \partial_x \mathcal{L}^0\mathcal{I}_{1/3}f(x,t) \\
&= 3\int_0^t A'\left( \frac{x}{(t-t')^{1/3}} \right) \frac{\mathcal{I}_{-1/3}f(t')}{(t-t')^{2/3}} \, dt'
\end{aligned}
\end{equation}
By Lemma \ref{L:Dbf2}, if $f\in C_0^\infty(\mathbb{R}^+)$, then $\mathcal{L}^{-1}f(x,t)$ is continuous in $x$ for all $x\in \mathbb{R}$ and, since $A'(0)=-(3\Gamma(\frac{1}{3}))^{-1}$, the second representation of $\mathcal{L}^{-1}f(x,t)$ in \eqref{E:Dbf2} gives
\begin{equation}
\label{E:Dbfeq4}
\mathcal{L}^{-1}f(0,t) = -f(t).
\end{equation}
By \eqref{E:Dbfeq1}, $\mathcal{L}^{-1}$ satisfies 
$$
\left\{
\begin{aligned}
& (\partial_t + \partial_x^3) \mathcal{L}^{-1}f(x,t) = 3\delta_0'(x)\mathcal{I}_{-1/3}f(t) && \text{for }(x,t)\in \mathbb{R}\times \mathbb{R}\\
& \mathcal{L}^{-1}f(x,0) = 0 && \text{for }x\in \mathbb{R}
\end{aligned}
\right.
$$
By Lemma \ref{L:Dbf2}, $\partial_x \mathcal{L}^0f(x,t)$ is continuous in $x$ for all $x\in \mathbb{R}$ and since $A'(0)=-(3\Gamma(\frac{1}{3}))^{-1}$,
\begin{equation}
\label{E:Dbfeq10}
\partial_x \mathcal{L}^0f(0,t) = -\mathcal{I}_{-1/3}f(t).
\end{equation}
Again by Lemma \ref{L:Dbf2},  $\partial_x \mathcal{L}^{-1}f(x,t)=\partial_x^2 \mathcal{L}^0\mathcal{I}_{1/3}f(x,t)$ is continuous in $x$ for $x\neq 0$ and has a step discontinuity of size $3\mathcal{I}_{-1/3}f(t)$ at $x=0$.  Since
\begin{align*}
\lim_{x\downarrow 0} \partial_x^2 \mathcal{L}^0f(x,t) &= -\int_0^{+\infty} \partial_y^3 \mathcal{L}^0f(y,t) \, dy \\
&= +\int_0^{+\infty} \mathcal{L}^0(\partial_t f)(y,t) \, dy && \text{by \eqref{AE:2006}} \\
&= 3 \int_{y=0}^{+\infty} A(y) \, dy \int_0^t  \partial_t\mathcal{I}_{-2/3}f(t') \, dt' && \text{by \eqref{E:Dbf} and Fubini}\\
&= \mathcal{I}_{-2/3}f(t)
\end{align*}
we have
\begin{equation}
\label{E:Dbfeq3}
\lim_{x\uparrow 0} \partial_x \mathcal{L}^{-1}f(x,t) = -2\mathcal{I}_{-1/3}f(t), \qquad \lim_{x\downarrow 0} \partial_x \mathcal{L}^{-1}f(x,t) = \mathcal{I}_{-1/3}f(t).
\end{equation}
By \eqref{E:Dbfeq2}, \eqref{E:Dbfeq4}, \eqref{E:Dbfeq10}, \eqref{E:Dbfeq3}, for yet to be assigned $h_1$ and $h_2$, we have
\begin{align}
\mathcal{L}^0h_1(0,t)+\mathcal{L}^{-1}h_2(0,t) &= h_1(t) - h_2(t) \label{E:Dbfeq5}\\
\lim_{x\uparrow 0} \mathcal{I}_{1/3}\partial_x (\mathcal{L}^0h_1(x,-)+\mathcal{L}^{-1}h_2(x,-))(t) &= -h_1(t) -2h_2(t) \label{E:Dbfeq6}\\
\lim_{x\downarrow 0} \mathcal{I}_{1/3}\partial_x (\mathcal{L}^0h_1(x,-)+\mathcal{L}^{-1}h_2(x,-))(t) &= -h_1(t) +h_2(t) \label{E:Dbfeq7}
\end{align}
If we are given $g_1(t)$, $g_2(t)$, $\phi$, and set
$$
\begin{bmatrix}
h_1 \\ h_2 
\end{bmatrix}
= \frac{1}{3}\begin{bmatrix}
2 & -1 \\
-1 & -1
\end{bmatrix}
\begin{bmatrix}
g_1 -e^{\cdot \partial_x^3}\phi\big|_{x=0} \\ 
\mathcal{I}_{1/3}(g_2-\partial_x e^{-\cdot \partial_x^3}\phi\big|_{x=0})
\end{bmatrix}
$$
then by letting $u(x,t)= e^{-t\partial_x^3}\phi(x) + \mathcal{L}^0h_1(x,t) + \mathcal{L}^{-1}h_2(x,t)$, we have 
$$
\left\{
\begin{aligned}
&(\partial_t + \partial_x^3)u(x,t) = 0 && \text{for }x\neq 0\\
&u(x,0) = \phi(x) && \text{for }x\in \mathbb{R}\\
&u(0,t) = g_1(t) && \text{for }t\in \mathbb{R}\\
&\lim_{x\uparrow 0}\partial_xu(x,t) = g_2(t) && \text{for }t\in \mathbb{R}\\
\end{aligned}
\right.
$$

Owing to the degeneracy in the right-hand limits \eqref{E:Dbfeq5}, \eqref{E:Dbfeq7}, we see that we cannot specify both boundary data $u(0,t)$ and derivative boundary data $\lim_{x\downarrow 0}\partial_x u(x,t)$ for the right half-line problem, which is consistent with the uniqueness calculation \eqref{E:1}.

\subsection{Nonlinear versions}

We define the \textit{Duhamel inhomogeneous solution operator} $\mathcal{D}$ as
\begin{equation}
\label{E:Di}
\mathcal{D}w(x,t) = \int_0^t e^{-(t-t')\partial_x^3}w(x,t') \, dt'
\end{equation}
so that 
\begin{equation}
\label{E:Dieq}
\left\{
\begin{aligned}
& (\partial_t + \partial_x^3) \mathcal{D}w(x,t) = w(x,t) && \text{for }(x,t) \in \mathbb{R}\times \mathbb{R} \\
& \mathcal{D}w(x,0) = 0 && \text{for }x\in \mathbb{R}
\end{aligned}
\right.
\end{equation}

For the right half-line problem \eqref{SE:100}, let
\begin{equation}
\label{E:370}
\Lambda_+ w = e^{-t\partial_x^3} \phi - \tfrac{1}{2}\mathcal{D}(\partial_x w^2) + \mathcal{L}^0h
\end{equation}
where
$$h(t) = f(t) - e^{-t\partial_x^3}\phi\big|_{x=0} + \tfrac{1}{2}\mathcal{D}(\partial_x w^2)(0,t)$$
and observe that if $u$ is such that $\Lambda_+ u = u$, then $u$ solves \eqref{SE:100}.
For the left half-line problem \eqref{SE:101}, let
\begin{equation}
\label{E:371}
\Lambda_- w = e^{-t\partial_x^3} \phi - \tfrac{1}{2}\mathcal{D}(\partial_x w^2) + \mathcal{L}^0h_1 + \mathcal{L}^{-1}h_2
\end{equation}
where
$$
\begin{bmatrix}
h_1(t) \\ h_2(t)
\end{bmatrix}
= \begin{bmatrix}
2 & 1 \\
-1 & -1
\end{bmatrix}
\begin{bmatrix}
g_1(t) - e^{-t\partial_x^3}\big|_{x=0} + \tfrac{1}{2}\mathcal{D}(\partial_xw^2)(0,t) \\ \mathcal{I}_{1/3}(g_2(\cdot) - \partial_x e^{-\cdot\partial_x^3}\phi\big|_{x=0} + \tfrac{1}{2}\partial_x \mathcal{D}(\partial_x w^2)(0,\cdot))
\end{bmatrix}
$$
and observe that if $u$ is such that $\Lambda_- u=u$, then $u$ solves \eqref{SE:101}.
One approach, then, to solving \eqref{SE:100} and \eqref{SE:101} is to prove that $\Lambda_+$, $\Lambda_-$ (or actually time-truncated versions of them) are contraction mappings in suitable Banach spaces.  As is the case for the IVP, we need the auxiliary Bourgain space \eqref{E:400}.

\begin{remark} In order to prove Lemma \ref{L:Dbf}\eqref{I:DbfBse}, we shall need to take $b<\frac{1}{2}$.  The $D_\alpha$ norm is a low frequency correction for the $X_{s,b}$ norm that is needed in order for the bilinear estimates (Lemma \ref{L:bilinear}) to hold for $b<\frac{1}{2}$.  This problem is particular to our treatment of initial-boundary value problems and does not arise in the standard treatment of the initial-value problem (IVP) using the $X_{s,b}$ spaces (see \cite{KPV96}).  In treating the IVP, one does not need the Duhamel boundary forcing operators and is thus at liberty to take $b>\frac{1}{2}$, and the bilinear estimate Lemma \ref{L:bilinear} holds in this case without the low frequency modification $D_\alpha$.
\end{remark}

Consider the space $Z$ consisting of all $w$ such that $w\in C(\mathbb{R}_t; H_x^s) \cap C(\mathbb{R}_x; H_t^\frac{s+1}{3}) \cap X_{s,b}\cap D_\alpha$ and $\partial_x w\in C(\mathbb{R}_x; H_t^\frac{s}{3})$.  Suppose we wanted to show that the maps $\Lambda_\pm$ above are contractions in a ball in $Z$ with radius determined by the norms of the initial and boundary data.   (This was done by \cite{CK02} for $\Lambda_+$ with $s=0$  without the estimates on $\partial_x u$ in $C(\mathbb{R}_x; H_t^\frac{s}{3})$, and their arguments easily extend to $-\frac{1}{2}<s<\frac{1}{2}$.)  The needed estimates for such an argument appear below in \S\ref{S:estimates} as Lemma \ref{L:G} for $e^{-t\partial_x^3}$, Lemma \ref{L:Di} for $\mathcal{D}$, Lemma \ref{L:Dbf} with $\lambda=0$ for $\mathcal{L}^0$, and Lemma \ref{L:Dbf} with $\lambda=-1$ for $\mathcal{L}^{-1}$.  The constraints in Lemma \ref{L:Dbf}\eqref{I:DbfBse} for $\lambda=0$ are $-\frac{1}{2}<s\leq 1$, and the constraints in Lemma \ref{L:Dbf}\eqref{I:DbfBse} for $\lambda=-1$ are $-\frac{3}{2}<s\leq 0$, thus restricting us to $-\frac{1}{2}<s\leq 0$.  In order to achieve the results in the wider range $-\frac{3}{4}<s<\frac{3}{2}$, $s\neq \frac{1}{2}$, we next introduce (in \S \ref{S:Dbf}) two analytic families of operators $\mathcal{L}_+^\lambda$ and $\mathcal{L}_-^\lambda$ such that $\mathcal{L}_\pm^0=\mathcal{L}^0$, $\mathcal{L}_\pm^{-1}=\mathcal{L}^{-1}$.  The solution properties are: 
$$
\left\{
\begin{aligned}
&(\partial_t + \partial_x^3)\mathcal{L}_+^\lambda f(x,t) = 3\frac{x_-^{\lambda-1}}{\Gamma(\lambda)} \mathcal{I}_{-\frac{2}{3}-\frac{\lambda}{3}}f(t) \\
&\mathcal{L}_+^\lambda f(x,0) = 0 \\
&\mathcal{L}_+^\lambda f(0,t) = e^{\pi i \lambda} f(t)
\end{aligned}
\right.
$$
and
$$
\left\{
\begin{aligned}
&(\partial_t + \partial_x^3)\mathcal{L}_-^\lambda f(x,t) = 3\frac{x_+^{\lambda-1}}{\Gamma(\lambda)} \mathcal{I}_{-\frac{2}{3}-\frac{\lambda}{3}}f(t) \\
&\mathcal{L}_-^\lambda f(x,0) = 0 \\
&\mathcal{L}_-^\lambda f(0,t) = 2\sin(\tfrac{\pi}{3}\lambda+\tfrac{\pi}{6}) f(t)
\end{aligned}
\right.
$$
Due to the support properties of $\frac{x_-^{\lambda-1}}{\Gamma(\lambda)}$ and $\frac{x_+^{\lambda-1}}{\Gamma(\lambda)}$, $(\partial_t + \partial_x^3)\mathcal{L}_+^{\lambda}f(x,t) = 0$ for $x>0$ and $(\partial_t + \partial_x^3)\mathcal{L}_-^\lambda f(x,t)=0$ for $x<0$.  For any $-\frac{3}{4}<s<\frac{3}{2}$, $s\neq \frac{1}{2}$, we will be able to address the right half-line problem \eqref{SE:100} by replacing $\mathcal{L}^0$ in \eqref{E:370} with $\mathcal{L}_+^\lambda$ for suitable $\lambda=\lambda(s)$ and address the left half-line problem \eqref{SE:101} by replacing $\mathcal{L}^0$, $\mathcal{L}^{-1}$ in \eqref{E:371} with $\mathcal{L}_-^{\lambda_1}$, $\mathcal{L}_-^{\lambda_2}$ for suitable $\lambda_1\neq \lambda_2$ chosen in terms of $s$.

After the classes $\mathcal{L}_\pm^{\lambda}$ have been defined and examined in \S \ref{S:Dbf}, some properties of the half-line Sobolev spaces $H_0^s(\mathbb{R}^+)$, $H^s(\mathbb{R}^+)$ will be given in \S\ref{S:N}.  The needed estimates for the contraction arguments are given in \S\ref{S:estimates}.  Finally in \S\ref{S:left}-\ref{S:linesegment}, we prove the local well-posedness results in Theorem \ref{T:main}.

\section{The Duhamel boundary forcing operator class}
\label{S:Dbf}

Define, for $\text{Re }\lambda >0$, and $f\in C_0^\infty(\mathbb{R}^+)$
\begin{equation} \label{AE:2020}
\begin{aligned}
\mathcal{L}_-^\lambda f(x,t) &= \left[ \frac{x_+^{\lambda-1}}{\Gamma(\lambda)} \ast \mathcal{L}^0(\mathcal{I}_{-\lambda/3}f)(-,t) \right](x) \\
&= \frac{1}{\Gamma(\lambda)} \int_{-\infty}^x (x-y)^{\lambda-1} \mathcal{L}^0(\mathcal{I}_{-\lambda/3}f)(y,t) \, dy
\end{aligned}
\end{equation}
and, with $\frac{x_-^{\lambda-1}}{\Gamma(\lambda)} = e^{i\pi\lambda}\frac{(-x)_+^{\lambda-1}}{\Gamma(\lambda)}$, define
\begin{equation} \label{AE:2020A}
\begin{aligned}
\mathcal{L}_+^\lambda f(x,t) &=  \left[ \frac{x_-^{\lambda-1}}{\Gamma(\lambda)} \ast \mathcal{L}^0(\mathcal{I}_{-\lambda/3}f)(-,t) \right](x) \\
&= \frac{e^{i\pi\lambda}}{\Gamma(\lambda)}  \int_x^{+\infty} (y-x)^{\lambda-1}\mathcal{L}^0(\mathcal{I}_{-\lambda/3}f)(y,t) \, dy 
\end{aligned}
\end{equation}
By integration by parts in \eqref{AE:2020}, the decay bounds provided by Lemma \ref{L:Dbf2}, and  \eqref{AE:2006},
\begin{align}
\mathcal{L}_-^\lambda f(x,t) &= \left[ \frac{x_+^{(\lambda+3)-1}}{\Gamma(\lambda+3)} \ast \partial_x^3\mathcal{L}^0f(-,t) \right](x) \notag \\
&= 
\begin{aligned}[t]
&3\frac{x_+^{(\lambda+3)-1}}{\Gamma(\lambda+3)}\mathcal{I}_{-\frac{2}{3}-\frac{\lambda}{3}}f(t) \\
&- \int_{-\infty}^x \frac{(x-y)^{(\lambda+3)-1}}{\Gamma(\lambda+3)} \mathcal{L}^0(\partial_t \mathcal{I}_{-\frac{\lambda}{3}}f)(y,t) \, dy 
\end{aligned}
\label{AE:2004}
\end{align}
For $\text{Re }\lambda >-3$, we may thus take \eqref{AE:2004} as the definition for $\mathcal{L}_-^\lambda f$. By integration by parts in \eqref{AE:2020A}, the decay bounds provided by Lemma \ref{L:Dbf2}, and  \eqref{AE:2006},
\begin{align}
\mathcal{L}_+^\lambda f(x,t) &=  \left[ \frac{x_-^{(\lambda+3)-1}}{\Gamma(\lambda+3)} \ast \partial_x^3\mathcal{L}f(-,t) \right](x) \notag \\
&= 
\begin{aligned}[t]
&3\frac{x_-^{(\lambda+3)-1}}{\Gamma(\lambda+3)}\mathcal{I}_{-\frac{2}{3}-\frac{\lambda}{3}}f(t) \\
&+ e^{i\pi\lambda} \int_{-\infty}^x \frac{(-x+y)^{(\lambda+3)-1}}{\Gamma(\lambda+3)} \mathcal{L}^0(\partial_t \mathcal{I}_{-\frac{\lambda}{3}}f)(y,t) \, dy
\end{aligned}
\label{AE:2004A}
\end{align}
For $\text{Re }\lambda >-3$, we may thus take \eqref{AE:2004A} as the definition for $\mathcal{L}_+^\lambda f$.
It is staightforward from these definitions that, in the sense of distributions
$$(\partial_t + \partial_x^3)\mathcal{L}_-^\lambda f(x,t) = 3\frac{x_+^{\lambda-1}}{\Gamma(\lambda)}\mathcal{I}_{-\frac{2}{3}-\frac{\lambda}{3}}f(t)$$
and
$$(\partial_t + \partial_x^3)\mathcal{L}_+^\lambda f(x,t) = 3\frac{x_-^{\lambda-1}}{\Gamma(\lambda)}\mathcal{I}_{-\frac{2}{3}-\frac{\lambda}{3}}f(t)$$

\begin{lemma}[Spatial continuity and decay properties for $\mathcal{L}_\pm^\lambda f(x,t)$]
\label{L:decayest}
Let $f\in C_0^\infty(\mathbb{R}^+)$, and fix $t\geq 0$.  We have 
$$\mathcal{L}_\pm^{-2}f=\partial_x^2\mathcal{L}^0\mathcal{I}_{\frac{2}{3}}f, \qquad \mathcal{L}_\pm^{-1}f=\partial_x\mathcal{L}^0\mathcal{I}_{\frac{1}{3}}f, \qquad \mathcal{L}_\pm^0f=\mathcal{L}f$$ 
Also, $\mathcal{L}_\pm^{-2}f(x,t)$ has a step discontinuity of size $3f(t)$ at $x=0$, otherwise for $x\neq 0$, $\mathcal{L}_\pm^{-2}f(x,t)$ is continuous in $x$.  For $\lambda >-2$, $\mathcal{L}_\pm^\lambda f(x,t)$ is continuous in $x$ for all $x\in \mathbb{R}$.  For $-2\leq \lambda \leq 1$, $0\leq t \leq 1$, $\mathcal{L}_-^\lambda f(x,t)$ satisfies the decay bounds
\begin{align*}
|\mathcal{L}_-^\lambda f(x,t)| &\leq c_{k,\lambda,f} \la x \ra^{-k}  && \forall \; x \leq 0, \quad \forall \; k\geq 0 \\
|\mathcal{L}_-^\lambda f(x,t)| &\leq c_{\lambda,f}\la x \ra^{\lambda-1} && \forall \; x \geq 0
\end{align*}
For $-2\leq \lambda \leq 1$, $0\leq t \leq 1$, $\mathcal{L}_+^\lambda f(x,t)$ satisfies the decay bounds
\begin{align*}
|\mathcal{L}_+^\lambda f(x,t)| &\leq c_{k,\lambda,f} \la x \ra^{-k}  && \forall \; x \geq 0, \quad \forall \; k\geq 0 \\
|\mathcal{L}_+^\lambda f(x,t)| &\leq c_{\lambda,f}\la x \ra^{\lambda-1} && \forall \; x \leq 0
\end{align*}
\end{lemma}
\begin{proof}
We only prove the bounds for $\mathcal{L}^\lambda_- f$, since the corresponding results for $\mathcal{L}^\lambda_+ f$ are obtained similarly.  For $x\leq -2$, the result follows by direct estimation in \eqref{AE:2004} using $|\mathcal{L}^0(\partial_t\mathcal{I}_{-\frac{\lambda}{3}}f)(y,t)| \leq c_{k,f}\la y \ra^{-k}\la x \ra^{-k}$ obtained from \eqref{AE:463} (since $|y|\geq |x|$).  Assume $x\geq 2$.  Let $\psi\in C^\infty(\mathbb{R})$ be such that $\psi(y)=1$ for $y\leq \frac{1}{4}$ and $\psi(y)=0$ for $y\geq \frac{3}{4}$.  Then
\begin{align*}
\mathcal{L}_-^\lambda f(x,t) &= \frac{x_+^{(\lambda+3)-1}}{\Gamma(\lambda+3)} \ast \partial_x^3 \mathcal{L}^0\mathcal{I}_{-\frac{\lambda}{3}}f(-,t) \\
&= 
\begin{aligned}[t]
&\int_{-\infty}^x \frac{(x-y)^{\lambda+2}}{\Gamma(\lambda+3)} \psi\left(\frac{y}{x}\right) \partial_y^3\mathcal{L}f^0\mathcal{I}_{-\frac{\lambda}{3}}(y,t) \, dy \\
&+ \int_{-\infty}^x \frac{(x-y)^{\lambda+2}}{\Gamma(\lambda+3)} \left[ 1 - \psi\left(\frac{y}{x}\right)\right] \partial_y^3\mathcal{L}^0\mathcal{I}_{-\frac{\lambda}{3}}f(y,t)  \, dy 
\end{aligned} \\
&= \text{I} + \text{II}
\end{align*}
In I, $y\leq \frac{3}{4}x$, integrate by parts,
\begin{align*}
\text{I} &= -\int_{-\infty}^x \partial_y^3 \left[ \frac{(x-y)^{\lambda+2}}{\Gamma(\lambda+3)} \psi\left(\frac{y}{x}\right)\right] \mathcal{L}^0\mathcal{I}_{-\frac{\lambda}{3}}f(y,t) \, dy \\
&= 
\begin{aligned}[t]
&\int_{-\infty}^x \frac{(x-y)^{\lambda-1}}{\Gamma(\lambda)} \psi\left( \frac{y}{x}\right) \mathcal{L}^0\mathcal{I}_{-\frac{\lambda}{3}}f(y,t) \, dy \\
&+\sum_{j=1}^3 c_j\int_{-\infty}^x \frac{(x-y)^{\lambda+j-1}}{\Gamma(\lambda+j)} \frac{1}{x^j} \psi^{(j)}\left( \frac{y}{x} \right) \mathcal{L}^0\mathcal{I}_{-\frac{\lambda}{3}}f(y,t) \, dy
\end{aligned}
\end{align*}
In the first of these terms, since $y\leq \frac{3}{4}x$, $(x-y)^{\lambda-1} \leq (\frac{1}{4})^{\lambda-1} x^{\lambda-1}$.  In the second term, $\frac{1}{4}x \leq y \leq \frac{3}{4}x$, and thus we can use the decay of $\mathcal{L}^0 \mathcal{I}_{-\lambda/3} f(y,t)$.
In II, $y\geq \frac{1}{4}x$, apply \eqref{AE:2006},
\begin{align*}
\text{II} &= \int_{-\infty}^x \frac{(x-y)^{\lambda+2}}{\Gamma(\lambda+3)} \left[ 1 - \psi\left(\frac{y}{x}\right)\right] ( 3\delta_0(y)\mathcal{I}_{-2/3}f(t)-\mathcal{L}^0(\partial_t\mathcal{I}_{-2/3}f)(y,t) ) \, dy \\
&= -\int_{-\infty}^x \frac{(x-y)^{\lambda+2}}{\Gamma(\lambda+3)} \left[ 1 - \psi\left(\frac{y}{x}\right)\right] \mathcal{L}^0(\partial_t\mathcal{I}_{-2/3}f)(y,t) \, dy
\end{align*}
Since $y\geq \tfrac{1}{4}x$, we have by Lemma \ref{L:Dbf2}, $$|\mathcal{L}^0(\partial_t \mathcal{I}_{-2/3}f)(y,t)| \leq c_k  \|f\|_{H^{2k+1}} \la x \ra^{-k} \la y \ra^{-k},$$ which establishes the bound.
\end{proof}

\begin{lemma}[Values of $\mathcal{L}_\pm^\lambda f(x,t)$ at $x=0$] 
\label{L:valueatzero}
For $\text{Re }\lambda>-2$,
\begin{equation}
\label{AE:2016}
\mathcal{L}_-^\lambda f(0,t) = 2\sin(\tfrac{\pi}{3}\lambda+\tfrac{\pi}{6})f(t)
\end{equation}
\begin{equation}
\label{AE:2016A}
\mathcal{L}_+^\lambda f(0,t) = e^{i\pi \lambda}f(t)
\end{equation}
\end{lemma}

In order to prove this, we need to compute the Mellin transform of each side of the Airy function.
\begin{lemma}[Mellin transform of the Airy function] \label{L:leftMellin}
If $0< \textnormal{Re }\lambda < \frac{1}{4}$, then
\begin{equation} \label{AE:2013}
\int_0^{+\infty} x^{\lambda-1}A(-x) \, dx = \tfrac{1}{3\pi}\Gamma(\lambda)\Gamma(-\tfrac{1}{3}\lambda + \tfrac{1}{3})\cos(\tfrac{2\pi}{3}\lambda - \tfrac{\pi}{6})
\end{equation}
 If $\textnormal{Re }\lambda >0 $, then
\begin{equation} \label{AE:2410}
\int_0^{+\infty} x^{\lambda-1} A(x) \, dx = \tfrac{1}{3\pi} \Gamma(\lambda) \Gamma(\tfrac{1}{3}-\tfrac{1}{3}\lambda) \cos( \tfrac{\pi}{3}\lambda+\tfrac{\pi}{6})
\end{equation}
Note that although $\Gamma(\frac{1}{3}-\frac{1}{3}\lambda)$ has poles at $\lambda=1,4,7, \cdots$, $\cos( \tfrac{\pi}{3}\lambda+\tfrac{\pi}{6})$ vanishes at these positions.
\end{lemma}
\begin{proof}
We shall only carry out the computation leading to \eqref{AE:2013}, since the one for \eqref{AE:2410} is similar.  Owing to the decay of the Airy function $A(-x) \leq c\la x \ra^{-1/4}$ for $x\geq 0$, the given expression is defined as an absolutely convergent integral.  In the calculation, we assume that $\lambda$ is real and $0<\lambda <\frac{1}{4}$, and by analyticity, this suffices to establish \eqref{AE:2013}.  Let $A_1(x) = \frac{1}{2\pi} \int_0^{+\infty} e^{ix\xi} e^{i\xi^3} \, d\xi$, so that $A(x)=2\text{Re }A_1(x)$.  Let $A_{1,\epsilon}(x) = \frac{1}{2\pi} \int_{\xi=0}^{+\infty} e^{ix\xi} e^{i\xi^3} e^{-\epsilon \xi} \, d\xi$.  Then, by dominated convergence and Fubini
\begin{align}
\indentalign \int_0^{+\infty} x^{\lambda-1}A_1(-x)\, dx \\
&= \lim_{\epsilon \downarrow 0} \lim_{\delta \downarrow 0} \int_{x=0}^{+\infty} x^{\lambda-1} e^{-\delta x} A_{1,\epsilon}(-x) \, dx  \notag\\
&=  \lim_{\epsilon \downarrow 0} \lim_{\delta \downarrow 0} \frac{1}{2\pi} \int_{\xi=0}^{+\infty} e^{i\xi^3}e^{-\epsilon\xi} \int_{x=0}^{+\infty} x^{\lambda-1}e^{-\delta x}e^{-ix\xi} \, dx \, d\xi . \label{E:301}
\end{align}
By a change of contour, 
\begin{equation}
\label{E:300}
\int_{x=0}^{+\infty} x^{\lambda-1}e^{-\delta x} e^{-ix\xi} \, dx = \xi^{-\lambda} e^{-\lambda \frac{\pi}{2}} \Gamma(\lambda, \delta/\xi)
\end{equation}
where $\Gamma(\lambda,z)= \int_{r=0}^{+\infty} r^{\lambda-1} e^{irz}e^{-r} \, dr$.  By dominated convergence, 
$$\lim_{\delta \downarrow 0} \int_{x=0}^{+\infty} x^{\lambda-1}e^{-\delta x} e^{-ix\xi} \, dx = \xi^{-\lambda} e^{-\lambda \frac{\pi}{2}} \Gamma(\lambda)$$
Since \eqref{E:300} is bounded independently of $\delta>0$, we have by dominated convergence
$$\eqref{E:301} =  \frac{1}{2\pi} \Gamma(\lambda) e^{-i\lambda \frac{\pi}{2}} \lim_{\epsilon \downarrow 0} \int_{\xi=0}^{+\infty} e^{i\xi^3}e^{-\epsilon\xi} \xi^{-\lambda} \, d\xi$$
Change variable $\eta=\xi^3$ and change contour, this becomes
$$ \frac{1}{6\pi} \Gamma(\lambda) e^{-\frac{2\pi \lambda i}{3}} e^{\frac{\pi i}{6}} \lim_{\epsilon \downarrow 0} \int_0^{+\infty} e^{-r} e^{-\epsilon( \frac{\sqrt{3}}{2}+i\frac{1}{2})r^{1/3}} r^{-\frac{2}{3}-\frac{\lambda}{3}} \, dr$$
Finally, dominated convergence yields
$$\int_0^{+\infty} x^{\lambda-1}A_1(-x)\, dx = \tfrac{1}{6\pi} e^{-\frac{2\pi \lambda i}{3}} e^{\frac{\pi i}{6}} \Gamma(\lambda)\Gamma(\tfrac{1}{3}-\tfrac{\lambda}{3})$$
  Using $A(x) = 2\text{Re }A_1(x)$, we obtain \eqref{AE:2013}
\end{proof}

Now we return to the proof of Lemma \ref{L:valueatzero}.

\begin{proof}[Proof of Lemma \ref{L:valueatzero}]
From \eqref{AE:2004}, 
\begin{equation} \label{AE:2015}
\mathcal{L}_-^\lambda f(0,t) = \int_{-\infty}^0 \frac{(-y)^{\lambda+2}}{\Gamma(\lambda+3)} \mathcal{L}^0(\partial_t \mathcal{I}_{-\frac{\lambda}{3}}f)(y,t) \, dy
\end{equation}
and from \eqref{AE:2004A},
\begin{equation} \label{AE:2015A}
\mathcal{L}_+^\lambda f(0,t) = e^{i\pi\lambda}\int_0^{+\infty} \frac{y^{\lambda+2}}{\Gamma(\lambda+3)} \mathcal{L}^0(\partial_t \mathcal{I}_{-\frac{\lambda}{3}}f)(y,t) \, dy
\end{equation}
By complex differentiation under the integral sign, \eqref{AE:2015} demonstrates that $\mathcal{L}_-^\lambda f(0,t)$ is analytic in $\lambda$ for $\text{Re }\lambda>-2$.  We shall only compute \eqref{AE:2016} for $0<\lambda<\tfrac{1}{4}$, $\lambda$ real.  By analyticity, the result will extend to the full range $\text{Re }\lambda>-2$.  For the computation in the range $0<\lambda<\tfrac{1}{4}$, we use the representation \eqref{AE:2020} in place of \eqref{AE:2015} to give
$$\mathcal{L}_-^\lambda f(0,t) = \int_{y=-\infty}^0 \frac{(-y)^{\lambda-1}}{\Gamma(\lambda)} \mathcal{L}^0f(y,t) \, dy$$
By  the decay for $A(-y)$, $y\geq 0$, we can apply Fubini to the above equation after inserting \eqref{E:Dbf} and then apply \eqref{AE:2013} to obtain
\begin{equation*} 
\mathcal{L}_-^\lambda f(0,t) = \tfrac{1}{\pi} \Gamma\left( -\tfrac{1}{3}\lambda+\tfrac{1}{3} \right) \Gamma\left( \tfrac{1}{3}\lambda + \tfrac{2}{3} \right) \cos\left( \tfrac{2\pi}{3}\lambda - \tfrac{\pi}{6} \right) \mathcal{I}_{\frac{1}{3}\lambda+\frac{2}{3}}(\mathcal{I}_{-\frac{\lambda}{3}-\frac{2}{3}}f)(t)
\end{equation*}
Using the identities $\Gamma(z)\Gamma(1-z)=\dfrac{\pi}{\sin \pi z}$, $\cos x = \sin( \frac{\pi}{2} - x)$, and $\sin 2x = 2\cos x \sin x$,
\begin{align*}
\mathcal{L}_-^\lambda f(0,t) &= \frac{\cos \left( \frac{2\pi}{3} \lambda -\frac{\pi}{6} \right)}{\sin \left( -\frac{\pi}{3}\lambda + \frac{\pi}{3} \right)} \mathcal{I}_{\frac{1}{3}\lambda+\frac{2}{3}}(h)(t) \\
&= 2 \sin \left( \tfrac{\pi}{3}\lambda + \tfrac{\pi}{6} \right)  \mathcal{I}_{\frac{1}{3}\lambda+\frac{2}{3}}(\mathcal{I}_{-\frac{\lambda}{3}-\frac{2}{3}}f)(t)
\end{align*}
giving \eqref{AE:2016}.
By complex differentiation under the integral sign, \eqref{AE:2015A} demonstrates that $f_+(t,\lambda)$ is analytic in $\lambda$ for $\text{Re }\lambda>-3$.  We shall only compute \eqref{AE:2016A} for $0<\lambda$, $\lambda$ real.  By analyticity, the result will extend to the full range $\text{Re }\lambda>-3$.  For the computation in the range $0<\lambda$, we use the representation \eqref{AE:2020A} in place of \eqref{AE:2015A} to give
$$\mathcal{L}_+^\lambda f(0,t) = e^{i\pi\lambda}\int_{y=0}^{+\infty} \frac{y^{\lambda-1}}{\Gamma(\lambda)} \mathcal{L}^0\mathcal{I}_{-\frac{\lambda}{3}}f(y,t) \, dy$$
By the decay of $A(y)$, $y\geq 0$, we can apply Fubini to obtain
$$\mathcal{L}_+^\lambda f(0,t) = \tfrac{1}{3\pi} \Gamma(\tfrac{1}{3}-\tfrac{1}{3}\lambda) \cos(\tfrac{\pi}{3}\lambda+ \tfrac{\pi}{6}) e^{i\pi\lambda} \mathcal{I}_{\frac{1}{3}\lambda+\frac{2}{3}}(\mathcal{I}_{-\frac{1}{3}\lambda-\frac{2}{3}}f)(t)$$
Using the same identities as above, we obtain \eqref{AE:2016A}.
\end{proof}

\section{Notations and some function space properties} \label{S:N}

  We use the notation $H^s$ to mean $H^s(\mathbb{R})$ (and not $H^s(\mathbb{R}^+)$ or $H_0^s(\mathbb{R}^+)$). The trace operator $\phi \mapsto\phi(0)$ is defined for $\phi\in H^s(\mathbb{R})$ when $s>\frac{1}{2}$.  For $s\geq 0$, define $\phi \in H^s(\mathbb{R}^+)$ if $\exists \; \tilde{\phi} \in H^s(\mathbb{R})$ such that $\tilde{\phi}(x)=\phi(x)$ for $x>0$; in this case we set $\|\phi\|_{H^s(\mathbb{R}^+)} = \inf_{\tilde{\phi}} \|\tilde{\phi} \|_{H^s(\mathbb{R})}$.  For $s \in \mathbb{R}$, define $\phi \in H_0^s(\mathbb{R}^+)$ if, when $\phi(x)$ is extended to $\tilde{\phi}(x)$ on $\mathbb{R}$ by setting $\tilde{\phi}(x)=0$ for $x<0$, then $\tilde{\phi}\in H^s(\mathbb{R})$; in this case we set $\| \phi \|_{H_0^s(\mathbb{R}^+)} = \|\tilde{\phi}\|_{H^s(\mathbb{R})}$.  For $s<0$, define $H^s(\mathbb{R}^+)$ as the dual space to $H_0^{-s}(\mathbb{R}^+)$, and define $H^s_0(\mathbb{R}^+)$ as the dual space to $H^{-s}(\mathbb{R}^+)$.  A definition for $H^s(0,L)$ can be given analogous to that for $H^s(\mathbb{R}^+)$. 

Define $\phi \in C_0^\infty(\mathbb{R}^+)$ if $\phi \in C^\infty(\mathbb{R})$ with $\text{supp }\phi \subset [0,+\infty)$ (so that, in particular, $\phi$ and all of its derivatives vanish at $0$), and $C_{0,c}^{\infty}(\mathbb{R}^+)$ as those members of $C_0^\infty(\mathbb{R}^+)$ with compact support.  We remark that $C_{0,c}^{\infty}(\mathbb{R}^+)$ is dense in $H_0^s(\mathbb{R}^+)$ for all $s\in \mathbb{R}$.

\begin{lemma}[\cite{CK02} Lemma 2.8]\label{CK28}  If $0\leq \alpha < \frac{1}{2}$, then $\| \theta h \|_{H^\alpha} \leq c\|h\|_{\dot{H}^\alpha}$ and $\| \theta h \|_{\dot{H}^{-\alpha}} \leq c\|h\|_{H^{-\alpha}}$, where $c=c(\alpha, \theta)$.
\end{lemma}

\begin{lemma}[\cite{JK95} Lemma 3.5] \label{JK35}
If $-\frac{1}{2}< \alpha< \frac{1}{2}$, then $\| \chi_{(0,+\infty)}f \|_{H^\alpha} \leq c \| f \|_{H^\alpha}$, where $c=c(\alpha)$.
\end{lemma}

\begin{lemma}[\cite{CK02} Prop.\ 2.4, \cite{JK95} Lemma 3.7, 3.8] \label{JK37}
If $\frac{1}{2}<\alpha<\frac{3}{2}$, then 
$$H_0^\alpha(\mathbb{R}^+)=\{ f\in H^\alpha(\mathbb{R}^+) \; \mid \; f(0)=0 \}.$$
If $\frac{1}{2}<\alpha<\frac{3}{2}$ and $f\in H^\alpha(\mathbb{R}^+)$ with $f(0)=0$, then $\|\chi_{(0,+\infty)} f\|_{H_0^\alpha(\mathbb{R}^+)} \leq c\|f\|_{H^\alpha(\mathbb{R}^+)}$, where $c=c(\alpha)$.
\end{lemma}

\begin{lemma}[\cite{CK02}, Lemma 5.1]
If $s\in \mathbb{R}$ and $0<b<1$, $0<\alpha<1$ then
$$\|\theta(t)w(x,t)\|_{X_{s,b}\cap D_\alpha} \leq c\|w\|_{X_{s,b}}$$
where $c=c(\theta)$.
\end{lemma}

\begin{lemma}[\cite{CK02} Cor.\ 2.1, Prop.\ 2.2] 
For $\alpha\geq 0$, $H_0^{-\alpha}(\mathbb{R}^+)$ is a complex interpolation scale.  For $\alpha\geq 0$, $H_0^\alpha(\mathbb{R}^+)$ is a complex interpolation scale.

\end{lemma}

\section{Estimates}
\label{S:estimates}

\subsection{Estimates for the Riemann-Liouville fractional integral}

In this section, we shall use the notation $\mathcal{J}_\alpha f =  \frac{t_+^{\alpha-1}}{\Gamma(\alpha)} \ast f $ for $f\in C_0^\infty(\mathbb{R})$ (no restriction on support of $f$ to $[0,+\infty)$).  This is in distinction to the definition of $\mathcal{I}_\alpha$, where we are convolving with a function $f$ supported in $[0,+\infty)$.

\begin{lemma}\label{L:RL5} Let $\alpha\in\mathbb{C}$.  If $\mu_1\in C_0^\infty(\mathbb{R})$ and $\mu_2\in C^\infty(\mathbb{R})$ such that $\mu_2=1$ on a neighborhood of $(-\infty, b]$, where $b=\sup \{ \, t \, | \, t\in \textnormal{supp }\mu_1 \, \}$, then $\mu_1 \mathcal{J}_\alpha \mu_2 h = \mu_1 \mathcal{J}_\alpha h$.  If $\mu_2\in C^\infty_0(\mathbb{R})$ and $\mu_1\in C^\infty(\mathbb{R})$ such that $\mu_1=1$ on a neighborhood of $[a, +\infty)$, where $a=\inf \{ \, t \, | \, t\in \textnormal{supp }\mu_2 \, \}$, then $\mu_1 \mathcal{J}_\alpha \mu_2 h = \mathcal{J}_\alpha \mu_2 h$
\end{lemma}
\begin{proof}
The first identity is clear from the integral definition if $\text{Re } \alpha>0$.  If $\text{Re } \alpha<0$, let $k\in \mathbb{N}$ be such that $-k<\text{Re }\alpha \leq -k+1$ so that $\mathcal{J}_\alpha=\partial_t^k\mathcal{J}_{\alpha+k}$.  Let $U$ be an open set such that 
$$\text{supp }\mu_1 \subset (-\infty, b] \subset U \subset \{ \, t\, | \, \mu_2(t)=1 \, \}$$
Then $\forall \; t\in U$,  $\mathcal{J}_{\alpha+k}h = \mathcal{J}_{\alpha+k}\mu_2 h$, which implies that $\forall \; t\in (-\infty,b]$, $\partial_t^k \mathcal{J}_{\alpha+k}h = \partial_t^k \mathcal{J}_{\alpha+k} \mu_2 h$, which implies that $\forall \; t\in\mathbb{R}$, $\mu_1\partial_t^k \mathcal{J}_{\alpha+k}h = \mu_1\partial_t^k \mathcal{J}_{\alpha+k} \mu_2 h$.  The second claim is clear by the integral definition if $\text{Re }\alpha >0$.  If $\text{Re }\alpha<0$, let $k\in \mathbb{N}$ be such that $-k<\text{Re }\alpha \leq -k+1$ so that $\mathcal{J}_\alpha=\mathcal{J}_{\alpha+k}\partial_t^k$.   Since $\text{supp }\partial_t^j \mu_2 \subset [a,+\infty)\subset \{ \, t\, | \, \mu_1(t)=1 \, \}$, we have 
$$\mu_1 \mathcal{J}_{\alpha+k} (\partial_t^j\mu_2) (\partial_t^{k-j}h) = \mathcal{J}_{\alpha+k}(\partial_t^j\mu_2)(\partial_t^{k-j}h)$$
and thus $\mu_1\mathcal{J}_{\alpha+k}\partial_t^k \mu_2 h = \mathcal{J}_{\alpha+k}\partial_t^k \mu_2 h$.
\end{proof}

\begin{lemma} \label{L:RL6} For $\gamma\in\mathbb{R}$, $s\in \mathbb{R}$, $\| \mathcal{J}_{i\gamma} h \|_{H^s(\mathbb{R})} \leq \cosh(\tfrac{1}{2}\pi\gamma)\| h \|_{H^s(\mathbb{R})}$
\end{lemma}
\begin{proof}
From \eqref{E:400}, we have 
$$\left( \frac{x_+^{i\gamma-1}}{\Gamma(i\gamma)} \right)\sphat(\xi) = 
\left\{
\begin{aligned}
&e^{\frac{1}{2}\pi \gamma} e^{-i\gamma\ln |\xi|} && \text{if }\xi >0 \\
&e^{-\frac{1}{2}\pi \gamma} e^{-i\gamma\ln |\xi|} && \text{if }\xi<0
\end{aligned}
\right.
$$
and thus $\left| \left( \frac{x_+^{i\gamma-1}}{\Gamma(i\gamma)} \right)\sphat(\xi) \right| \leq 2\cosh(\frac{1}{2}\pi \gamma)$.
\end{proof}

\begin{lemma} \label{L:RL1}
If $0\leq \textnormal{Re }\alpha < +\infty$ and $s\in \mathbb{R}$, then
\begin{gather} 
\| \mathcal{I}_{-\alpha}h \|_{H_0^s(\mathbb{R}^+)}\leq c e^{\frac{1}{2}\textnormal{Im }\alpha}\|h\|_{H_0^{s+\alpha}(\mathbb{R}^+)} \label{AE:315} \\
\| \mathcal{J}_{-\alpha}h \|_{H^s(\mathbb{R})}\leq c e^{\frac{1}{2}\textnormal{Im }\alpha} \|h\|_{H^{s+\alpha}(\mathbb{R})} \label{AE:316}
\end{gather}
\end{lemma}
\begin{proof}
\eqref{AE:316} is immediate from \eqref{E:410}.  \eqref{AE:315} then follows from \eqref{AE:316} by Lemma \ref{L:RL} and a density argument.
\end{proof}

\begin{lemma} \label{L:RL2}
If $0\leq \textnormal{Re }\alpha < +\infty$, $s\in \mathbb{R}$, $\mu, \mu_2 \in C_0^\infty(\mathbb{R})$
\begin{align}
\| \mu \mathcal{I}_\alpha h \|_{H_0^s(\mathbb{R}^+)} &\leq c e^{\frac{1}{2}\textnormal{Im } \alpha}\|h\|_{H_0^{s-\alpha}(\mathbb{R}^+)} & c&=c(\mu)\label{AE:318} \\
\| \mu \mathcal{J}_\alpha \mu_2 h \|_{H^s(\mathbb{R})} &\leq c e^{\frac{1}{2}\textnormal{Im } \alpha} \|h\|_{H^{s-\alpha}(\mathbb{R})}  &  c&= c(\mu,\mu_2)\label{AE:319}
\end{align}
where $c=c(\mu, \mu_2)$.
\end{lemma}
\begin{proof}  We first explain how \eqref{AE:318} follows from \eqref{AE:319}. Given $\mu$, let $b= \sup\{ \, t \, | \, t \in \text{supp } \mu \, \}$.  Take $\mu_2\in C_0^\infty(\mathbb{R})$, $\mu_2=1$ on $[0,b]$. Then, when restricting to $h\in C_0^\infty(\mathbb{R}^+)$, we have $\mu \mathcal{I}_\alpha h = \mu \mathcal{J}_\alpha \mu_2 h$.  By Lemma \ref{L:RL} and a density argument, we obtain \eqref{AE:318}.  Now we prove \eqref{AE:319}.  We first need the special case $s=0$.\\
\textit{Claim}.  If  $k\in \mathbb{Z}_{\geq 0}$, then $\|\mu \mathcal{J}_k \mu_2 h\|_{L^2(\mathbb{R})} \leq c \| h \|_{H^{-k}(\mathbb{R})}$, where $c=c(\mu, \mu_2)$.\\
To prove this claim, consider $k\in \mathbb{N}$.  If $g\in C_0^\infty(\mathbb{R})$ with $\| g\|_{L^2}\leq 1$, then
\begin{align*}
\| \mu \mathcal{J}_k \mu_2 h \|_{L^2} &= \frac{1}{\Gamma(k)}\sup_g \int_t \mu(t) \int_{s=-\infty}^t(t-s)^{k-1} \mu_2(s)h(s) \, ds \, g(t) \, dt \\
&= \frac{1}{\Gamma(k)} \sup_g \int_s h(s) \, \mu_2(s) \int_{t=s}^{+\infty} \mu(t)(t-s)^{k-1}g(t)\, dt \, ds \\
&\leq \frac{1}{\Gamma(k)}\| h \|_{H^{-k}} \left\| \mu_2(s) \int_{t=s}^{+\infty} \mu(t)(t-s)^{k-1}g(t)\, dt  \right\|_{H^k(ds)} \\
& \leq c\| h\|_{H^{-k}} \|g\|_{L^2}
\end{align*}
The case $k =0$ is trivial, concluding the proof of the claim.\\
To prove \eqref{AE:319}, we first take $\alpha=k\in \mathbb{Z}_{\geq 0}$, $s=m\in\mathbb{Z}$, $h\in C_0^\infty(\mathbb{R})$.\\
\textit{Case 1}.  $m\geq 0$.
\begin{align*}
\|\mu \mathcal{J}_k \mu_2 h \|_{H^m} &\leq \| \mu \mathcal{J}_k \mu_2 h \|_{L^2} + \sum_{j=0}^m \| \mu^{(j)} \mathcal{J}_{k-m+j} \mu_2 h \|_{L^2} \\
&\leq c(\| h\|_{H^{-k}} + \sum_{j=0}^m\|h\|_{H^{m-k-j}})\leq c\|h\|_{H^{m-k}}
\end{align*}
by appealing to the claim or Lemma \ref{L:RL1}.\\
\textit{Case 2}.  $m<0$.
Let $\mu_3=1$ on $\text{supp }\mu$, $\mu_3\in C_0^\infty(\mathbb{R}^+)$.  
$$\mu \mathcal{J}_k \mu_2 h = \mu \partial_t^{-m} \mathcal{J}_{k-m} \mu_2 h = \mu \partial_t^{-m} \mu_3 \mathcal{J}_{k-m} \mu_2 h$$
and therefore
$$\| \mu \mathcal{J}_k \mu_2 h \|_{H^m} \leq \| \mu_3 \mathcal{J}_{k-m} \mu_2 h \|_{L^2}$$
and we conclude by applying the claim.  

Next, we extend to $\alpha=k+i\gamma$ for $k,\gamma \in \mathbb{R}$, as follows.  Let $\mu_3=1$ on a neighborhood of $(-\infty,b]$, where $b=\sup \{ \, t \, | \, t\in \text{supp }\mu \, \}$, and let $\mu_4=1$ on a neighborhood of $[a,+\infty)$, where $a=\inf \{ \, t \, | \, t\in \text{supp }\mu_2 \, \}$, so that $\mu_3\mu_4 \in C_0^\infty(\mathbb{R})$.  By Lemma \ref{L:RL5},
$$\mu \mathcal{J}_{k+i\gamma} \mu_2 h = \mu \mathcal{J}_{i\gamma} \mu_3 \mu_4 \mathcal{J}_k \mu_2 h$$
By Lemma \ref{L:RL6},
$$\| \mu \mathcal{J}_{k+i\gamma} \mu_2 h \|_{H^m} \leq c \cosh (\tfrac{1}{2}\pi \gamma)\| \mu_3 \mu_4 \mathcal{J}_k \mu_2 h\|_{H^m}$$
which is bounded as above.  We can now apply interpolation to complete the proof.
\end{proof}

\subsection{Estimates for the group}

The operator $e^{-t\partial_x^3}$ was defined above in \eqref{E:G} satisfying \eqref{E:Geq}.

\begin{lemma} \label{L:G}
Let $s\in \mathbb{R}$.  Then
\begin{enumerate}
\item \label{I:Gst} \textnormal{(Space traces)} $\| e^{-t\partial_x^3}\phi(x)\|_{C(\mathbb{R}_t; H^s_x)} \leq c\| \phi\|_{H^s}$.
\item \label{I:Gtt} \textnormal{(Time traces)} $\| \theta(t) e^{-t\partial_x^3} \phi(x) \|_{C(\mathbb{R}_x; H_t^\frac{s+1}{3})} \leq c\|\phi\|_{H^s}$.
\item \label{I:Gdtt} \textnormal{(Derivative time traces)} $\| \theta(t) \partial_x e^{-t\partial_x^3} \phi(x) \|_{C(\mathbb{R}_x; H_t^\frac{s}{3})} \leq c\|\phi\|_{H^s}$.
\item \label{I:GBse} \textnormal{(Bourgain space estimate)} If $0<b<1$ and $0<\alpha<1$, then $\| \theta(t) e^{-t\partial_x^3} \phi(x) \|_{X_{s,b}\cap D_\alpha} \leq c\|\theta\|_{H^1}\|\phi\|_{H^s}$, where $c$ is independent of $\theta$.
\end{enumerate}
\end{lemma}
\begin{proof}
\eqref{I:Gst},\eqref{I:GBse} follow from the definition \eqref{E:G} and \eqref{I:Gtt},\eqref{I:Gdtt} appear in \cite{KPV91}.
\end{proof}

\subsection{Estimates for the Duhamel inhomogeneous solution operator}

The operator $\mathcal{D}$ was defined above in \eqref{E:Di} satisfying \eqref{E:Dieq}.

Let
$$\|u\|_{Y_{s,b}} = \left( \iint_{\xi,\tau} \la\tau\ra^{2s/3}\la\tau-\xi^3\ra^{2b} |\hat{u}(\xi,\tau)|^2 \, d\xi \, d\tau \right)^{1/2}$$

\begin{lemma} \label{L:Di}
Let $s\in \mathbb{R}$.  Then
\begin{enumerate}
\item \label{I:Dist} \textnormal{(Space traces)} If $0\leq b < \frac{1}{2}$, then $$ \| \theta(t) \mathcal{D}w(x,t) \|_{C(\mathbb{R}_t; H^s_x)} \leq c\|w\|_{X_{s,-b}}.$$
\item \label{I:Ditt} \textnormal{(Time traces)} If $0<b<\frac{1}{2}$, then 
\begin{align*}
\indentalign 
\| \theta(t) \mathcal{D}w(x,t) \|_{C(\mathbb{R}_x; H^\frac{s+1}{3}_t)} \\
&\leq \left\{
\begin{aligned}
& c\|w\|_{X_{s,-b}} && \textnormal{if }-1\leq s\leq \tfrac{1}{2} \\
& c(\|w\|_{X_{s,-b}}+ \|w\|_{Y_{s,-b}}) && \textnormal{for any }s
\end{aligned}
\right.
\end{align*}
If $s<\frac{7}{2}$, then $\|\theta(t)\mathcal{D}w(x,t) \|_{C(\mathbb{R}_x;H_0^\frac{s+1}{3}(\mathbb{R}_t^+))}$ has the same bound.
\item \label{I:Didtt} \textnormal{(Derivative time traces)} If $0<b<\frac{1}{2}$, then 
\begin{align*}
\indentalign 
\| \theta(t) \partial_x \mathcal{D}w(x,t) \|_{C(\mathbb{R}_x; H^\frac{s}{3}_t)} \\
&\leq \left\{
\begin{aligned}
& c\|w\|_{X_{s,-b}} && \textnormal{if }0\leq s\leq \tfrac{3}{2} \\
& c(\|w\|_{X_{s,-b}}+ \|w\|_{Y_{s,-b}}) && \textnormal{for any }s
\end{aligned}
\right.
\end{align*} 
If $s<\frac{9}{2}$, then $\|\theta(t)\partial_x\mathcal{D}w(x,t) \|_{C(\mathbb{R}_x;H_0^\frac{s}{3}(\mathbb{R}_t^+))}$ has the same bound.
\item \label{I:DiBse} \textnormal{(Bourgain space estimate)} If $0\leq b < \frac{1}{2}$ and $\alpha\leq 1-b$, then $\ds \|\theta(t) \mathcal{D}w(x,t) \|_{X_{s,b}\cap D_\alpha} \leq c\|w\|_{X_{s,-b}}$.
\end{enumerate}
\end{lemma}

\begin{remark} The need for the $Y_{s,b}$ (time-adapted) Bourgain space arises here in Lemma \ref{L:Di}\eqref{I:Ditt}\eqref{I:Didtt} in order the cover the full interval $-\frac{3}{4}<s<\frac{3}{2}$ ($s\neq \frac{1}{2}$).  It is, however, only an intermediate device since the bilinear estimate in Lemma \ref{L:bilinear}\eqref{I:BYX} enables us to avoid carrying out the contraction argument in $Y_{s,b}$.
\end{remark}

\begin{proof}
\eqref{I:DiBse} is Lemma 5.4 in \cite{CK02} (although $Y_b$ has a different definition from ours) and \eqref{I:Dist} is a standard estimate (see the techniques of Lemmas 5.4, 5.5 in \cite{CK02}).  \eqref{I:Ditt} is Lemma 5.5 in \cite{CK02} and the proof of \eqref{I:Didtt} is modelled on the proof of Lemma 5.5 in \cite{CK02}.
\end{proof}

\subsection{Estimates for the Duhamel boundary forcing operator class}

The operators $\mathcal{L}_\pm^\lambda$ were defined above in \eqref{E:Dbf} solving \eqref{E:Dbfeq1}, \eqref{E:Dbfeq2}.

\begin{lemma} \label{L:Dbf}
Let $s\in \mathbb{R}$.  Then
\begin{enumerate}
\item \label{I:Dbfst} \textnormal{(Space traces)} If $s-\frac{5}{2}<\lambda<s+\frac{1}{2}$, $\lambda<\frac{1}{2}$, and $\textnormal{supp }f\subset[0,1]$, then $\| \mathcal{L}_\pm^\lambda f(x,t) \|_{C(\mathbb{R}_t; H^s_x)} \leq c\|f\|_{H_0^\frac{s+1}{3}(\mathbb{R}^+)}$.
\item \label{I:Dbftt} \textnormal{(Time traces)} If $-2< \lambda <1$, then 
$$\|\theta(t) \mathcal{L}_\pm^\lambda f(x,t) \|_{C(\mathbb{R}_x; H_0^\frac{s+1}{3}(\mathbb{R}_t^+))} \leq c\|f\|_{H_0^\frac{s+1}{3}(\mathbb{R}^+)}.$$
\item \label{I:Dbfdtt} \textnormal{(Derivative time traces)} If $-1 < \lambda < 2$, then  $$\| \theta(t) \partial_x \mathcal{L}_\pm^\lambda f (x,t) \|_{C(\mathbb{R}_x; H_0^\frac{s}{3}(\mathbb{R}_t^+))} \leq c\|f\|_{H_0^\frac{s+1}{3}(\mathbb{R}_t^+)}.$$
\item \label{I:DbfBse} \textnormal{(Bourgain space estimate)} If $s-1\leq \lambda < s+\frac{1}{2}$, $\lambda<\frac{1}{2}$, $\alpha\leq \frac{s-\lambda+2}{3}$, and $0\leq b<\frac{1}{2}$, then $\|\theta(t) \mathcal{L}_\pm^\lambda f(x,t) \|_{X_{s,b}\cap D_\alpha} \leq c\|f\|_{H_0^\frac{s+1}{3}(\mathbb{R}_t^+)}$.
\end{enumerate}
\end{lemma}

\begin{remark}  The restrictions on $s$, $\lambda$ in Lemma \ref{L:Dbf}\eqref{I:Dbfst}\eqref{I:DbfBse} are the primary purpose for introducing the analytic families $\mathcal{L}_\pm^\lambda$ and not simply using $\mathcal{L}^0$ for the right half-line problem and $\mathcal{L}^0$, $\mathcal{L}^{-1}$ for the left half-line problem. Note that by the assumption $\lambda<s+\frac{1}{2}$, we have $\frac{s-\lambda+2}{3}>\frac{1}{2}$, and thus we may take $\frac{1}{2}<\alpha\leq \frac{s-\lambda+2}{3}$, which is needed in order to meet the hypotheses of the bilinear estimates in Lemma \ref{L:bilinear}.
\end{remark}

\begin{proof}
We restrict to $\mathcal{L}_-^\lambda$ for notational convenience.  Also, we assume in the proof that $f\in C_0^\infty(\mathbb{R}^+)$.  The estimates, of course, extend by density.
To prove (a), we use ($\widehat{\quad}$ denoting the Fourier transform in $x$ alone)
$$( \mathcal{L}^\lambda f)\sphat(\xi,t) = (\xi-i0)^{-\lambda} \int_0^t e^{i(t-t')\xi^3} \mathcal{I}_{-\frac{\lambda}{3}-\frac{2}{3}}f(t') \, dt'$$
By the change of variable $\eta=\xi^3$ and the support properties of $ \mathcal{I}_{-\frac{\lambda}{3}-\frac{2}{3}}f(t')$, 
\begin{align*}
\|\phi\|_{H^s}^2 &\leq \int_\eta |\eta|^{-\frac{2\lambda}{3}-\frac{2}{3}} \la \eta \ra^{\frac{2s}{3}} \left| \int_0^t e^{i(t-t')\eta} \mathcal{I}_{-\frac{\lambda}{3}-\frac{2}{3}} f(t') \, dt' \right|^2 \, d\eta \\
&= \int_\eta |\eta|^{-\frac{2\lambda}{3}-\frac{2}{3}} \la \eta \ra^{\frac{2s}{3}} |( \chi_{(-\infty,t)}\mathcal{I}_{-\frac{\lambda}{3}-\frac{2}{3}}f)\sphat(\eta)|^2 \, d\eta
\end{align*}
noting that $\lambda <\frac{1}{2} \Longrightarrow -\frac{2}{3}\lambda-\frac{2}{3}$ and $s-\frac{5}{2}<\lambda<s+\frac{1}{2} \Longrightarrow -1<-\frac{2\lambda}{3}-\frac{2}{3}+\frac{2s}{3}<1$.
By Lemma \ref{CK28} (to replace $|\eta|^{-\frac{2\lambda}{3}-\frac{2}{3}}$ by $\la \eta \ra^{-\frac{2\lambda}{3}-\frac{2}{3}}$), Lemma \ref{JK35} (to remove the time cutoff factor $\chi_{(-\infty,t)}$), and Lemma \ref{L:RL1} (to estimate $\mathcal{I}_{-\frac{\lambda}{3}-\frac{2}{3}}$) we obtain the estimate in (a).

To prove (b), we first note that the change of variable $t' \rightarrow t-t'$ shows that
$$(I-\partial_t^2)^\frac{s+1}{6} \int_{-\infty}^t e^{-(t-t')\partial_x^3} h(t') \, dt' = \int_{-\infty}^t e^{-(t-t')\partial_x^3} (I-\partial_t^2)^{\frac{s+1}{6}} h(t') \, dt'$$
and thus (b) is equivalent to
$$\left\| \int_\xi e^{ix\xi}(\xi-i0)^{-\lambda} \int_{-\infty}^t e^{+i(t-t')\xi^3} (\mathcal{I}_{-\frac{\lambda}{3}-\frac{2}{3}}f)(t') \, dt' \, d\xi \right\|_{L_t^2} \leq c\|f\|_{L_t^2}$$
Using that $\chi_{(-\infty,t)} = \frac{1}{2}\sgn(t-t')+\frac{1}{2}$,
\begin{align*}
\indentalign 
\int_\xi e^{ix\xi}(\xi-i0)^{-\lambda} \int_{-\infty}^t e^{+i(t-t')\xi^3} (\mathcal{I}_{-\frac{\lambda}{3}-\frac{2}{3}}f)(t') \, dt' \, d\xi \\
&= 
\begin{aligned}[t]
&\int_\tau e^{it\tau} \left[ \lim_{\epsilon \downarrow 0} \int_{|\tau-\xi^3|> \epsilon} e^{ix\xi} \frac{(\tau-i0)^{\frac{\lambda}{3}+\frac{2}{3}}(\xi-i0)^{-\lambda}}{\tau-\xi^3} \, d\xi \right] \hat{f}(\tau) \, d\tau\\
&+\int_\xi e^{ix\xi}(\xi-i0)^{-\lambda} \int_{-\infty}^{+\infty} e^{+i(t-t')\xi^3} (\mathcal{I}_{-\frac{\lambda}{3}-\frac{2}{3}}f)(t') \, dt' \, d\xi
\end{aligned} \\
&= \text{I}+ \text{II}
\end{align*}
We can rewrite II as
$$\text{II} = \int_\xi e^{ix\xi} (\mathcal{I}_{-\frac{\lambda}{3}-\frac{2}{3}}f)\sphat(\xi^3) (\xi-i0)^{-\lambda} e^{it\xi^3} \, d\xi$$
The substitution $\eta =\xi^3$ and \eqref{E:410} gives
$$\text{II} = \int_\eta e^{it\eta} e^{ix\eta^{1/3}} (\eta-i0)^{\frac{\lambda}{3}+\frac{2}{3}} (\eta^{1/3}-i0)^{-\lambda} \eta^{-2/3}\hat{f}(\eta) \, d\eta$$
which is clearly $L^2_t \to L^2_t$ bounded.  In addressing term I, it suffices to show that
\begin{equation} \label{AE:434}
\lim_{\epsilon\downarrow 0} \int_{|\tau-\xi^3|>\epsilon} e^{ix\xi} \frac{(\tau-i0)^{\frac{\lambda}{3}+\frac{2}{3}}(\xi-i0)^{-\lambda}}{\tau-\xi^3} \, d\xi
\end{equation}
is bounded independently of $\tau$.  Changing variable $\xi \rightarrow \tau^{1/3}\xi$, and using that 
 $$(\tau^{1/3}\xi-i0)^{-\lambda} = \tau_+^{-\lambda/3}(c_1\xi_+^{-\lambda}+c_2\xi_-^{-\lambda}) + \tau_-^{-\lambda/3}(c_1\xi_-^{-\lambda}+c_2\xi_+^{-\lambda})$$
we get
$$\eqref{AE:434}=\chi_{\tau > 0} \int_\xi e^{i\tau^{1/3}x\xi} \frac{c_1\xi_+^{-\lambda}+c_2\xi_-^{-\lambda}}{1-\xi^3} \, d\xi + \chi_{\tau < 0} \int_\xi e^{i\tau^{1/3}x\xi} \frac{c_1 \xi_-^{-\lambda} + c_2 \xi_+^{-\lambda}}{1-\xi^3}$$
The treatment of both integrals is similar, so we will only consider the first of the two.  Let $\psi(\xi)=1$ near $\xi=1$, and $0$ outside $[\frac{1}{2},\frac{3}{2}]$.  Then this term breaks into 
$$c_1\int_\xi e^{ix\tau^{1/3}\xi} \frac{\psi(\xi)\xi_+^{-\lambda}}{1-\xi^3} \, d\xi + \int_\xi e^{ix\tau^{1/3}\xi} \frac{(1-\psi(\xi)) (c_1 \xi_+^{-\lambda} + c_2 \xi_-^{-\lambda})}{1-\xi^3} \, d\xi = \text{I}_a+\text{I}_b$$
The integrand in term $\text{I}_b$ is an $L^1$ function (provided $\lambda>-2$), so $|\text{Term I}_b| \leq c$.  Term $\text{I}_a$ is
$$c_1 \int_\xi e^{ix\tau^{1/3}\xi} \frac{\psi(\xi) \xi_+^{-\lambda}}{1+\xi+\xi^2} \frac{1}{1-\xi} \, d\xi$$
This becomes convolution of a Schwartz class function with a phase shifted $\text{sgn }x$ function, which is bounded on $L_t^2$, completing the proof of (b).  

Part (c) of the theorem is a corollary of (b) and the fact that $\partial_x \mathcal{L}_\pm^\lambda = \mathcal{L}_\pm^{\lambda-1} \mathcal{I}_{1/3}$.

To prove (d), first note that by \eqref{E:410}
$$( \mathcal{L}_-^\lambda f)\sphat (\xi,t) =  (\xi-i0)^{-\lambda} \int_\tau \frac{ e^{it\tau}-e^{it\xi^3}}{\tau-\xi^3} (\tau-i0)^{\frac{\lambda}{3}+\frac{2}{3}}\hat{f}(\tau) \, d\tau$$
Let $\psi(\tau)\in C^\infty(\mathbb{R})$ such that $\psi(\tau)=1$ for $|\tau|\leq 1$ and $\psi(\tau)=0$ for $|\tau|\geq 2$.  Set
$$\hat{u}_1(\xi,t) =  (\xi-i0)^{-\lambda} \int_\tau \frac{ e^{it\tau}-e^{it\xi^3}}{\tau-\xi^3} \psi(\tau-\xi^3)(\tau-i0)^{\frac{\lambda}{3}+\frac{2}{3}}\hat{f}(\tau) \, d\tau$$
$$\hat{u}_{2,1}(\xi,t) = (\xi-i0)^{-\lambda} \int_\tau \frac{ e^{it\tau}}{\tau-\xi^3} (1-\psi(\tau-\xi^3))(\tau-i0)^{\frac{\lambda}{3}+\frac{2}{3}}\hat{f}(\tau) \, d\tau$$
$$\hat{u}_{2,2}(\xi,t) =  (\xi-i0)^{-\lambda} \int_\tau \frac{ e^{it\xi^3}}{\tau-\xi^3} (1-\psi(\tau-\xi^3))(\tau-i0)^{\frac{\lambda}{3}+\frac{2}{3}}\hat{f}(\tau) \, d\tau$$
so that $\mathcal{L}_-^\lambda f = u_1 + u_{2,1} + u_{2,2}$.
For $-1<\lambda<\frac{1}{2}$, both $(\xi-i0)^{-\lambda}$, $(\tau-i0)^{\frac{\lambda}{3}+\frac{2}{3}}$ are square integrable functions and thus
\begin{equation}
\label{E:411}
\|u_{2,1}\|_{X_{s,b}}^2 \leq c \int_\tau |\tau|^{\frac{2\lambda}{3}+\frac{4}{3}} \left( \int_\xi \frac{|\xi|^{-2\lambda}\la \xi \ra^{2s}}{\la \tau-\xi^3 \ra^{2-2b}} \, d\xi \right) |\hat{f}(\tau)|^2 \, d\tau
\end{equation}
Since $-1<\lambda<\frac{1}{2}$, we have $-1<-\frac{2\lambda}{3}-\frac{2}{3}<0$ and
\begin{equation}
\label{E:412}
\int_\xi \frac{|\xi|^{-2\lambda}\la \xi \ra^{2s}}{\la \tau-\xi^3 \ra^{2-2b}} \, d\xi = \int_\eta |\eta|^{-\frac{2\lambda}{3}-\frac{2}{3}} \la \eta \ra^{\frac{2s}{3}} \la \tau-\eta \ra^{-2+2b} \, d\eta \leq c\la \tau \ra^{-\frac{2\lambda}{3}-\frac{2}{3}+\frac{2s}{3}}
\end{equation}
This is obtained by separately considering the cases $|\eta|\leq 1$, $|\tau|<<|\eta|$, and $|\eta|<<|\tau|$, and using that $s-1 \leq \lambda < s+\frac{1}{2}$ implies $-1<\frac{2s}{3}-\frac{2\lambda}{3}-\frac{2}{3}\leq 0$.  Combining \eqref{E:410} and \eqref{E:411} gives the appropriate bound for $\|u_{2,1}\|_{X_{s,b}}$.  
To address the term $u_{2,2}$, we first note that $u_{2,2}(x,t) = \theta(t)e^{-t\partial_x^3}\phi(x)$, where
\begin{equation}
\label{E:414}
\hat{\phi}(\xi) = (\xi-i0)^{-\lambda} \int_\tau \frac{1-\psi(\tau-\xi^3)}{\tau-\xi^3} (\mathcal{I}_{-\frac{\lambda}{3}-\frac{2}{3}}f)\sphat(\tau) \, d\tau
\end{equation}
Taking $h=\mathcal{I}_{-\frac{\lambda}{3}-\frac{2}{3}}f$ (so that $h\in C_0^\infty(\mathbb{R}^+)$ by Lemma \ref{L:RL}), we claim that
\begin{equation}
\label{E:413}
\int_\tau \hat{h}(\tau) \frac{1-\psi(\tau-\xi^3)}{\tau-\xi^3} \, d\tau = \int_\tau \hat{h}(\tau) \beta(\tau-\xi^3) \, d\tau
\end{equation}
where $\beta\in \mathcal{S}(\mathbb{R})$.  This follows from the fact that $\text{supp }h\subset [0,+\infty)$ as follows:  
 Let $ \hat{g_1}(\tau)=\frac{1-\psi(-\tau)}{\tau}$.  Then
$$ g_1(t)=\tfrac{i}{2}\text{sgn } t - \tfrac{i}{4\pi}\int_s \text{sgn}(t-s)\hat{\psi}(s) \, ds$$
Let $\alpha\in C^\infty(\mathbb{R})$ be such that $\alpha(t)=1$ for $t>0$ and $\alpha(t)=-1$ for $t<-1$, and set
$$g_2(t)=\tfrac{i}{2}\alpha(t) - \tfrac{i}{4\pi}\int_s \text{sgn}(t-s)\hat{\psi}(s) \, ds$$
 To show that $g_2\in \mathcal{S}(\mathbb{R})$, note that by the definition and the fact that $\hat{\psi}\in \mathcal{S}$, we have $g_2\in C^\infty(\mathbb{R})$.  If $t>0$, then since $\tfrac{1}{2\pi}\int\hat{\psi}(\tau)\, d\tau = \psi(0)=1$, we have 
$$
g_2(t)=\tfrac{i}{2} - \tfrac{i}{4\pi}\int_s \text{sgn}(t-s)\hat{\psi}(s) \, ds = \tfrac{i}{2\pi} \int_{s>t} \hat{\psi}(s) \, ds
$$
If $t<-1$, then likewise we have
$$
 g_2(t)=-\tfrac{i}{2} - \tfrac{i}{4\pi}\int_s \text{sgn}(t-s)\hat{\psi}(s) \, ds  = \tfrac{i}{2\pi} \int_{s<t} \hat{\psi}(s) \, ds$$
which provide the decay at $\infty$ estimates for $g_2$ and all of its derivatives, establishing that $g_2\in \mathcal{S}(\mathbb{R})$.  Since $g_1(t)=g_2(t)$ for $t>0$ and $h\in C_0^\infty(\mathbb{R}^+)$ we have
\begin{align*}
\indentalign \int_\tau \hat{h}(\tau) \frac{1-\psi(\tau-\xi^3)}{\tau-\xi^3} \, d\tau = -(\hat{h}*\hat{g_1})(\xi^3)  = -2\pi \widehat{hg_1}(\xi^3) \\
& = -2\pi \widehat{hg_2}(\xi^3)  = \int_\tau \hat{h}(\tau) \beta(\tau-\xi^3) \, d\tau
\end{align*}
where $\beta(\tau)=-\hat{g_2}(-\tau)$, and $\beta\in\mathcal{S}(\mathbb{R})$ since $g_2\in \mathcal{S}(\mathbb{R})$, thus establishing \eqref{E:413}. 
To complete the treatment of $u_{2,2}$, it suffices to show, by Lemma \ref{L:G}\eqref{I:GBse}, that $\|\phi\|_{H^s} \leq c\|f\|_{H^{\frac{s+1}{3}}}$.
By \eqref{E:414}, \eqref{E:413}, Cauchy-Schwarz and the fact that $|\beta(\tau-\xi^3)| \leq c\la \tau-\xi^3 \ra^{-N}$ for $N>>0$,
\begin{align*}
\|\phi\|_{H^s} &\leq \int_\xi \la \xi \ra^{2s} |\xi|^{-2\lambda} \left( \int_\tau \beta(\tau-\xi^3) |\tau|^{\frac{\lambda}{3}+\frac{2}{3}}|\hat{f}(\tau)| \, d\tau \right)^2 \, d\xi \\
&\leq \int_\tau \left( \int_\xi |\xi|^{-2\lambda} \la \xi \ra^{2s} \la \tau-\xi^3 \ra^{-2N+2} \, d\xi \right) |\tau|^{\frac{2\lambda}{3}+\frac{4}{3}} |\hat{f}(\tau)|^2 \, d\tau
\end{align*}
After the change of variable $\eta=\xi^3$, the inner integral becomes ($\lambda<\frac{1}{2} \Longrightarrow -\frac{2\lambda}{3}-\frac{2}{3}>-1$)
$$ \int_\eta |\eta|^{-\frac{2\lambda}{3}-\frac{2}{3}} \la \eta \ra^{\frac{2s}{3}} \la \tau -\eta \ra^{-2N+2} \, d\eta \leq c\la \tau \ra^{-\frac{2\lambda}{3}-\frac{2}{3}+\frac{2s}{3}}$$
This latter estimate can be obtained by considering cases $|\eta|\leq 1$, $|\eta|\leq \frac{1}{2}|\tau|$, and $|\eta| \geq \frac{1}{2}|\tau|$ (using $-1\leq \lambda-s \Longrightarrow -\frac{2\lambda}{3}-\frac{2}{3}+\frac{2s}{3}\leq 0$).  By the power series expansion for $e^{it(\tau-\xi^3)}$, $u_1(x,t) = \sum_{k=1}^{+\infty} \frac{1}{k!} \theta_k(t) e^{-t\partial_x^3} \phi_k(x)$, where $\theta_k(t) = i^kt^k\theta(t)$ and
$$\hat{\phi}_k(\xi) = (\xi-i0)^{-\lambda} \int_\tau (\tau-\xi^3)^{k-1} \psi(\tau-\xi^3) ( \mathcal{I}_{-\frac{2}{3}-\frac{\lambda}{3}} f)\sphat (\tau) \, d\tau$$
By Lemma \ref{L:G} \eqref{I:GBse}, it suffices to show that $\|\phi_k \|_{H^s} \leq c \|f\|_{H^\frac{s+1}{3}}$.  Note that
\begin{align*}
\|\phi_k\|_{H^s} &\leq \int_\xi \la \xi \ra^{2s} |\xi|^{-2\lambda} \left( \int_{|\tau-\xi^3| \leq 1} |\tau|^{\frac{\lambda}{3}+\frac{2}{3}} |\hat{f}(\tau)| \, d\tau \right)^2 \, d\xi \\
& \leq  \int_\tau \left( \int_{|\tau-\xi^3|\leq 1} \la \xi \ra^{2s} |\xi|^{-2\lambda} \, d\xi \right) |\tau|^{\frac{2\lambda}{3}+\frac{4}{3}} |\hat{f}(\tau)|^2 \, d\tau
\end{align*}
The substitution $\eta=\xi^3$ on the inner integral provides the needed bound.
\end{proof}

\subsection{Bilinear estimates}

\begin{lemma} \label{L:bilinear}
\begin{enumerate}
\item \label{I:BXX} For $s>-\frac{3}{4}$, $\exists \; b=b(s)<\frac{1}{2}$ such that $\forall \; \alpha>\frac{1}{2}$, we have 
\begin{equation}
\label{AE:230}
\|\partial_x(uv)\|_{X_{s,-b}} \leq c\|u\|_{X_{s,b}\cap D_{\alpha}}\|v\|_{X_{s,b}\cap D_\alpha}
\end{equation}
\item \label{I:BYX} For $-\frac{3}{4}<s<3$, $\exists \; b=b(s)<\frac{1}{2}$ such that $\forall \; \alpha>\frac{1}{2}$, we have 
\begin{equation}
\label{AE:3151}
\|\partial_x(uv)\|_{Y_{s,-b}} \leq c\|u\|_{X_{s,b}\cap D_\alpha}\|v\|_{X_{s,b}\cap D_\alpha}
\end{equation}
\end{enumerate}
\end{lemma}

\begin{remark}  The purpose of introducing the $D_\alpha$ low frequency correction factor is to validate the bilinear estimates above for $b<\frac{1}{2}$.  Recall that the need to take $b<\frac{1}{2}$ arose in Lemma \ref{L:Dbf}\eqref{I:DbfBse}.  
\end{remark}

We shall prove Lemma \ref{L:bilinear} by the calculus techniques of \cite{KPV96}.  We begin with some elementary integral estimates.

\begin{lemma} \label{calc1}
If $\frac{1}{4}<b<\frac{1}{2}$, then 
\begin{equation}
\int_{-\infty}^{+\infty} \frac{dx}{\la x-\alpha\ra^{2b}\la x-\beta\ra^{2b}} \leq \frac{c}{\la \alpha-\beta\ra^{4b-1}}
\end{equation}
\end{lemma}

\begin{proof}
By translation, it suffices to prove the inequality for $\beta=0$.  One then treats the cases $|\alpha|\leq 1$ and $|\alpha|\geq 1$ separately, and for the latter case, uses $\la x - \alpha \ra^{-2b} \la x \ra ^{-2b} \leq |x-\alpha|^{-2b}|x|^{-2b}$ and scaling.
\end{proof}

The following is \cite{KPV96} Lemma 2.3 (2.11) with $2b-\frac{1}{2}=1-l$ verbatim.
\begin{lemma} \label{calc2}
If $b<\frac{1}{2}$, then
\begin{equation}
\int_{|x|\leq \beta} \frac{dx}{\la x\ra^{4b-1}|\alpha-x|^{1/2}} \leq \frac{c(1+\beta)^{2-4b}}{\la \alpha\ra^\frac{1}{2}}
\end{equation} 
\end{lemma} 

\begin{proof}[Proof of Lemma \ref{L:bilinear} \eqref{I:BXX}]
We begin by addressing $-\frac{3}{4}<s<-\frac{1}{2}$.  The proof is modelled on the proof for $b>\frac{1}{2}$ given by \cite{KPV96}.  Essentially, we only need to replace one of the calculus estimates (\cite {KPV96} Lemma 2.3 (2.8)) in that paper with a suitable version for $b<\frac{1}{2}$ (Lemma \ref{calc1}).  Let $\rho=-s$.  It suffices to prove
\begin{equation} \label{AE:2100}
\iint_\ast \frac{|\xi| d(\xi,\tau)}{\la \tau-\xi^3 \ra^b \la \xi \ra^\rho} \frac{ \la \xi_1 \ra^\rho \hat{g}_1(\xi,\tau_1)}{\beta(\xi_1,\tau_1)} \frac{\la \xi_2 \ra^\rho\hat{g}_2(\xi_2,\tau_2)}{\beta(\xi_2,\tau_2)} \leq c \|d\|_{L^2} \|g_1\|_{L^2} \|g_2\|_{L^2}
\end{equation}
for $\hat{d}\geq 0$, $\hat{g_1}\geq 0$, $\hat{g}_2\geq 0$, where $\ast$ indicates integration over $\xi$, $\xi_1$, $\xi_2$, subject to the constraint $\xi=\xi_1+\xi_2$, and over $\tau$, $\tau_1$, $\tau_2$, subject to the constraint $\tau=\tau_1+\tau_2$, and where $\beta_j(\xi_j,\tau_j) = \la \tau_j-\xi_j^3\ra^b + \chi_{|\xi_j|\leq 1}\la \tau_j \ra^\alpha$.  By symmetry, it suffices to consider the case $|\tau_2-\xi_2^3| \leq |\tau_1-\xi_1^3|$.  We address \eqref{AE:2100} in pieces by the Cauchy-Schwarz method of \cite{KPV96}. We shall assume that $|\xi_1|\geq 1$ and $|\xi_2|\geq 1$, since otherwise, the bound \eqref{AE:2100} reduces to the case $\rho=0$, which has already been established in \cite{CK02}.\\

\noindent \textit{Case 1}.  If $|\tau_2-\xi_2^3|\leq |\tau_1-\xi_1^3| \leq |\tau-\xi^3|$, then we shall show
\begin{equation} \label{AE:2101}
\frac{|\xi|}{\la \tau-\xi^3 \ra^b \la \xi \ra^\rho} \left( \iint_{\tau_1,\xi_1} \frac{\la \xi_1 \ra^{2\rho} \la \xi_2 \ra^{2\rho}}{ \la \tau_1-\xi_1^3 \ra^{2b} \la \tau_2-\xi_2^3 \ra^{2b}} \, d\xi_1 \, d\tau_1 \right)^{1/2} \leq c
\end{equation}
To prove this, we note that 
\begin{equation} \label{AE:2102}
\tau-\xi^3 + 3\xi\xi_1\xi_2 = (\tau_2-\xi_2^3) + (\tau_1-\xi_1^3)
\end{equation}
By lemma \ref{calc1} with $\alpha =\xi_1^3$ and $\beta=\xi_1^3+\tau-\xi^3+3\xi\xi_1\xi_2$, we get that \eqref{AE:2101} is bounded by
$$
 \frac{|\xi|}{\la \tau-\xi^3 \ra^b \la \xi \ra^\rho} \left( \int_{\xi_1} \frac{\la \xi_1 \ra^{2\rho} \la \xi_2 \ra^{2\rho}}{\la \tau-\xi^3 + 3\xi\xi_1\xi_2 \ra^{4b-1}} \, d\xi_1 \right)^{1/2}
$$
By \eqref{AE:2102}, $|\xi\xi_1\xi_2|\leq |\tau-\xi^3|$.  Substituting $|\xi_1\xi_2| \leq |\tau-\xi^3| |\xi|^{-1}$ into the above gives that it is bounded by
\begin{equation} \label{AE:2104}
\frac{|\xi|^{1-\rho} \la \tau-\xi^3 \ra^{\rho-b}}{\la \xi \ra^\rho} \left( \int_{\xi_1} \frac{d\xi_1}{\la \tau-\xi^3+3\xi\xi_1\xi_2 \ra^{4b-1}} \right)^{1/2}
\end{equation}
Let $u=\tau-\xi^3 + 3\xi\xi_1\xi_2$, so that, by \eqref{AE:2102}, we have $|u| \leq 2|\tau-\xi^3|$.  The corresponding differential is
$$d\xi_1 = \frac{cdu}{|\xi|^{1/2} |u-(\tau-\frac{1}{4}\xi^3)|^{1/2}}$$
Substituting into \eqref{AE:2104}, we obtain that \eqref{AE:2104} is bounded by
$$
\frac{ |\xi|^{\frac{3}{4}-\rho} \la \tau-\xi^3 \ra^{\rho-b}}{\la \xi \ra^\rho} \left( \int_{|u|\leq 2|\tau-\xi^3|} \frac{du}{\la u \ra^{4b-1} |u-(\tau-\frac{1}{4}\xi^3)|^{1/2}} \right)^{1/2}
$$
By Lemma \ref{calc2}, this is controlled by
$$\frac{\la \tau-\xi^3 \ra^{\rho+1-3b}}{\la \xi \ra^{2\rho-\frac{3}{4}}\la \tau-\frac{1}{4}\xi^3 \ra^{1/4}}$$
This expression is bounded, provided $b\geq \frac{1}{9}\rho + \frac{5}{12}$.\\
\noindent \textit{Case 2}.  $|\tau_2-\xi_2^3| \leq |\tau_1-\xi_1^3|$, $|\tau-\xi^3|\leq |\tau_1-\xi_1^3|$.  In this case, we shall prove the bound
\begin{equation} \label{AE:2300}
\frac{1}{\la \tau_1-\xi_1^3 \ra^b} \left( \iint_{\xi,\tau} \frac{|\xi|^{2-2\rho} |\xi\xi_1\xi_2|^{2\rho}}{\la \xi \ra^{2\rho} \la \tau-\xi^3 \ra^{2b} \la \tau_2-\xi_2^3 \ra^{2b}} \, d\xi \, d\tau \right)^{1/2} \leq c
\end{equation}
Since 
\begin{equation} \label{AE:2301}
(\tau_1-\xi_1^3)+(\tau_2-\xi_2^3)-(\tau-\xi^3) = 3\xi\xi_1\xi_2
\end{equation}
we have, by Lemma \ref{calc1} with $\alpha=\xi^3$, $\beta=\xi^3+(\tau_1-\xi_1^3)-3\xi\xi_1\xi_2$, that \eqref{AE:2300} is bounded by
\begin{equation} \label{AE:2302}
\frac{1}{\la \tau_1-\xi_1^3 \ra^b} \left( \int_\xi \frac{ \la \xi \ra^{2-4\rho} |\xi\xi_1\xi_2|^{2\rho}}{\la \tau_1-\xi_1^3 -3\xi\xi_1\xi_2\ra^{4b-1}} \, d\xi \right)^{1/2}
\end{equation}
We address \eqref{AE:2302} in cases.  Cases 2A and 2B differ only in the bound used for $\la \xi \ra^{2-4\rho}$, while Case 2C is treated somewhat differently. \\
\noindent \textit{Case 2A}.  $|\xi_1|\sim |\xi|$ or $|\xi_1|<<|\xi|$.  Here, we use $\la \xi \ra^{2-4\rho} \leq \la \xi_1 \ra^{2-4\rho}$.\\
\noindent \textit{Case 2B}.  $|\xi|<<|\xi_1|$ and [$|\tau_1|>>\frac{1}{4}|\xi_1|^3$ or $|\tau_1|<<\frac{1}{4}|\xi_1|^3$].  Here, we use $\la \xi \ra^{2-4\rho} \leq 1$.\\
\noindent \textit{Cases 2A and 2B}.  In the setting of Case 2A, let $g(\xi_1)=\la \xi_1\ra^{1-2\rho}$, and in the setting of Case 2B, let $g(\xi_1) = 1$.  Since by \eqref{AE:2301}, $|\xi\xi_1\xi_2|\leq |\tau_1-\xi_1^3|$, \eqref{AE:2302} is bounded by
\begin{equation} \label{AE:2305}
g(\xi_1) \la \tau_1-\xi_1^3\ra^{\rho-b} \left( \int_\xi \frac{d\xi}{\la \tau_1-\xi_1^3-3\xi\xi_1\xi_2 \ra^{4b-1}} \right)^{1/2}
\end{equation}
Set $u=\tau_1-\xi_1^3-3\xi\xi_1\xi_2$.  Then 
$$du=3\xi_1(\xi_1-2\xi)d\xi=c|\xi_1|^{1/2}|u-(\tau_1-\tfrac{1}{4}\xi_1^3)|^{1/2}d\xi$$
which, upon substituting in \eqref{AE:2305}, gives that it is bounded by
$$ \frac{g(\xi_1)\la \tau_1-\xi_1^3 \ra^{\rho-b}}{|\xi|^{1/4}} \left( \int_{|u|\leq 2|\tau_1-\xi_1^3|} \frac{du}{\la u \ra^{4b-1}|u-(\tau_1-\frac{1}{4}\xi_1^3)|^{1/2}} \right)^{1/2}
$$
By Lemma \ref{calc2}, this is controlled by
\begin{equation} \label{AE:2307}
 \frac{g(\xi_1) \la \tau_1-\xi_1^3 \ra^{\rho+1-3b}}{|\xi_1|^{1/4}\la \tau_1-\frac{1}{4}\xi_1^3 \ra^{1/4}}
\end{equation}
In Case 2A, $g(\xi_1) = \la \xi_1 \ra^{1-2\rho}$, and \eqref{AE:2307} becomes
$$\frac{\la \tau_1-\xi_1^3 \ra^{\rho+1-3b}}{\la \xi_1 \ra^{2\rho-\frac{3}{4}} \la \tau_1-\frac{1}{4}\xi_1^3 \ra^{1/4}}$$
which is bounded provided $b>\frac{1}{9}\rho+\frac{5}{12}$.  In Case 2B, $g(\xi_1) = 1$, and \eqref{AE:2307} becomes
$$ \frac{ \la \tau_1-\xi_1^3 \ra^{\rho+1-3b}}{\la \xi_1 \ra^{1/4}\la \tau_1-\frac{1}{4}\xi_1^3 \ra^{1/4}}$$
which is bounded (under the restrictions of Case 2B) provided $b\geq \frac{1}{3}\rho + \frac{1}{4}$.\\
\noindent \textit{Case 2C}.  $|\xi|<<|\xi_1|$ and $|\tau_1|\sim \frac{1}{4}|\xi_1|^3$.  Here, we return to \eqref{AE:2302} and use that $|\tau_1|\sim \frac{1}{4}|\xi_1|^3$ and $3|\xi\xi_1\xi_2| \leq \frac{1}{4}|\xi_1|^3$ implies $\la \tau_1-\xi_1^3 -3\xi\xi_1\xi_2 \ra \sim \la \xi_1 \ra^3$.  Substituting into \eqref{AE:2302}, we find that it is bounded by
$$ \la \xi_1\ra^{3\rho-15b+3} \left( \int_{|\xi|\leq |\xi_1|} \la \xi \ra^{2-4\rho} \, d\xi \right)^{1/2} \leq \la \xi_1 \ra^{\rho-15b+\frac{9}{2}}$$
which is bounded provided $b\geq \frac{1}{15}\rho+\frac{3}{10}$.

We have completed the proof for $-\frac{3}{4}<s<-\frac{1}{2}$, and we shall now extend this result to all $s>-\frac{3}{4}$ by interpolation.  From the above, we have \eqref{AE:230} for $s=-\frac{5}{8}$ and some $b<\frac{1}{2}$.  As a consequence,
\begin{align*}
\indentalign \| \partial_x (uv) \|_{X_{\frac{3}{8},-b}} \\
&\leq \| \partial_x (uv) \|_{X_{-\frac{5}{8},-b}} + \| \partial_x [ (\partial_x u) v] \|_{X_{-\frac{5}{8},-b}} + \| \partial_x [ u (\partial_xv) ] \|_{X_{-\frac{5}{8},-b}} \\
& \leq (\| u \|_{X_{-\frac{5}{8},b}\cap D_\alpha} + \| \partial_x u \|_{X_{-\frac{5}{8},b}\cap D_\alpha}) (\| v \|_{X_{-\frac{5}{8},b}\cap D_\alpha} + \| \partial_x v \|_{X_{-\frac{5}{8},b}\cap D_\alpha})\\
& \leq \| u \|_{X_{\frac{3}{8},b}\cap D_\alpha} \| v \|_{X_{\frac{3}{8},b}\cap D_\alpha}
\end{align*}
thus establishing \eqref{AE:230} for $s=\frac{3}{8}$.  Now we can interpolate between the cases $s=-\frac{5}{8}$ and $s=\frac{3}{8}$ to obtain \eqref{AE:230} for $-\frac{3}{4} < s \leq \frac{3}{8}$.  Similarly, we can extend \eqref{AE:230} to all $s>-\frac{3}{4}$.  
\end{proof}

\begin{proof}[Proof of Lemma \ref{L:bilinear}\eqref{I:BYX}]
  First we address the range $-\frac{1}{2}<s<-\frac{3}{4}$. Let $\rho=-s$. Note that by the $X_{s,b}$ bilinear estimate Lemma \ref{L:bilinear}\eqref{I:BXX}, it suffices to prove the lemma under the assumption $|\tau| \leq \frac{1}{8}|\xi|^3$.  Constant multiples are routinely omitted from the calculation. \\

\noindent \textit{Step 1}.  If $|\xi_1|\geq 1$, $|\xi_2|\geq 1$, $|\tau_2-\xi_2^3|\leq |\tau_1-\xi_1^3|$, $|\tau_1-\xi_1^3|\leq 1000|\tau-\xi^3|$, and $|\tau|\leq \frac{1}{8}|\xi|^3$, then the expression
\begin{equation} \label{AE:200}
\frac{|\xi|}{\la \tau \ra^\frac{\rho}{3} \la \xi \ra^{3b}} \left( \int_{\xi_1} \int_{\tau_1} \frac{|\xi_1|^{2\rho} |\xi_2|^{2\rho}}{\la \tau_1-\xi_1^3 \ra^{2b} \la \tau_2-\xi_2^3 \ra^{2b}} \, d\tau_1 \, d\xi_1 \right)^{1/2}
\end{equation}
is bounded.

\noindent \textit{Proof.} Applying Lemma \ref{calc1}, using $\tau_2-\xi_2^3=(\tau-\xi^3)-(\tau_1-\xi_1^3)+3\xi\xi_1\xi_2$, we get that \eqref{AE:200} is bounded by
$$\frac{|\xi|}{\la \tau \ra^\frac{\rho}{3} \la \xi \ra^{3b}} \left( \int_{\xi_1}  \frac{|\xi_1|^{2\rho} |\xi_2|^{2\rho}}{\la \tau-\xi^3+3\xi\xi_1\xi_2 \ra^{4b-1}} \, d\xi_1 \right)^{1/2}
$$
Using that $|\xi_1||\xi_2| \leq \dfrac{|\tau-\xi^3|}{|\xi|}$, this is controlled by
\begin{equation} \label{AE:202}
 \frac{|\xi|^{1-\rho}|\tau-\xi^3|^\rho}{\la \xi \ra^{3b} \la \tau \ra^{\rho/3}} \left( \int_{\xi_1}  \frac{1}{\la \tau-\xi^3+3\xi\xi_1\xi_2 \ra^{4b-1}} \, d\xi_1 \right)^{1/2}
\end{equation}
Set 
$$
u=\tau-\xi^3+3\xi\xi_1(\xi-\xi_1)
$$
so that $3\xi(\xi_1-\frac{1}{2}\xi)^2=u-(\tau-\frac{1}{4}\xi^3)$, and thus
$$\tfrac{3}{\sqrt{2}}|\xi||2\xi_1-\xi| = |\xi|^{1/2}|u-(\tau-\tfrac{1}{4}\xi^3)|^{1/2}$$
Also, $du=3\xi(\xi-2\xi_1)\, d\xi_1$.  It follows from the hypotheses of this step that the range of integration is a subset of $|u| \leq |\tau-\xi^3|$.  With this substitution, we see that \eqref{AE:202} is bounded by
$$ \frac{|\xi|^{1-\rho}|\tau-\xi^3|^\rho}{\la \xi \ra^{3b} \la \tau \ra^{\rho/3}} \left( \int_{|u|\leq |\tau-\xi^3|} \frac{du}{\la u \ra^{4b-1} |\xi|^{1/2} |u-(\tau-\tfrac{1}{4}\xi^3)|^{1/2}} \right)^{1/2}
$$
By Lemma \ref{calc2}, this is controlled by
$$\frac{|\xi|^{\frac{3}{4}-\rho}|\tau-\xi^3|^\rho \la \tau-\xi^3 \ra^{1-2b}}{\la \xi \ra^{3b} \la \tau \ra^{\rho/3} \la \tau-\frac{1}{4}\xi^3 \ra^{1/4}}
$$
If $|\tau|\leq \frac{1}{8}|\xi|^3$, then this reduces to
$$\frac{|\xi|^{\frac{3}{4}-\rho} \la\xi\ra^{3\rho} \la \xi \ra^{3(1-2b)}}{\la \xi \ra^{3b} \la \xi \ra^{3/4}}$$
and the exponent $2\rho-9b+3\leq 0$ provided $b\geq \frac{2}{9}\rho+\frac{1}{3}$.\\

\noindent \textit{Step 2}.  If $|\xi_1|\geq 1$, $|\xi_2|\geq 1$, $|\tau_2-\xi_2^3| \leq |\tau_1-\xi_1^3|$, $|\tau-\xi^3| \leq \frac{1}{1000}|\tau_1-\xi_1^3|$, and $|\tau|\leq \frac{1}{8}|\xi|^3$, then
\begin{equation} \label{AE:204}
\frac{|\xi_1|^\rho}{\la \tau_1-\xi_1^3 \ra^b} \left( \int_\xi \int_\tau \frac{|\xi|^2 |\xi_2|^{2\rho}}{\la \tau \ra^{2\rho/3} \la \xi \ra^{6b} \la \tau_2-\xi_2^3 \ra^{2b}} \, d\xi \, d\tau \right)^{1/2}
\end{equation}
is bounded.  \\
\textit{Proof.}  Since $|\tau|\leq|\xi|^3$, we have $\dfrac{1}{\la \xi \ra^{6b-2\rho}} \leq \dfrac{1}{\la \tau \ra^{2b-\frac{2\rho}{3}}}$, and thus \eqref{AE:204} is bounded by
$$ \frac{|\xi_1|^\rho}{\la \tau_1-\xi_1^3 \ra^b} \left( \int_\xi \int_\tau \frac{|\xi|^2 |\xi_2|^{2\rho}}{\la \xi \ra^{2\rho} \la \tau \ra^{2b}  \la \tau_2-\xi_2^3 \ra^{2b}} \, d\xi \, d\tau \right)^{1/2}
$$
Carrying out the $\tau$ integral and applying Lemma \ref{calc1}, we see that this is controlled by
\begin{equation}\label{AE:210}
 \frac{|\xi_1|^\rho}{\la \tau_1-\xi_1^3 \ra^b} \left( \int_\xi  \frac{|\xi|^2 |\xi_2|^{2\rho}}{\la \xi \ra^{2\rho} \la \tau_1-\xi_1^3 -3\xi\xi_1\xi_2 + \xi^3 \ra^{4b-1}} \, d\xi  \right)^{1/2}
\end{equation}
\textit{Case 1.}  $3|\xi\xi_1\xi_2| \leq \frac{1}{2}|\tau_1-\xi_1^3|$. \\ Since $|\tau-\xi^3|<<|\tau_1-\xi_1^3|$ and $|\tau|\leq \frac{1}{8}|\xi|^3$, we have $|\xi|^3 << |\tau_1-\xi_1^3|$, giving $$\la \tau_1-\xi_1^3-3\xi\xi_1\xi_2 + \xi^3 \ra \sim \la \tau_1-\xi_1^3 \ra$$ and thus \eqref{AE:210} is bounded by
$$ \frac{|\xi_1|^\rho}{\la \tau_1-\xi_1^3 \ra^{3b-\frac{1}{2}}} \left( \int_\xi \frac{ |\xi|^2 |\xi_2|^{2\rho}}{\la \xi \ra^{2\rho}} \, d\xi \right)^{1/2}
$$
Using that $|\xi\xi_1\xi_2| \leq |\tau_1-\xi_1^3|$, this is controlled by
\begin{equation} \label{AE:207}
\frac{|\tau_1-\xi_1^3|^\rho}{\la \tau_1-\xi_1^3 \ra^{3b-\frac{1}{2}}} \left( \int_\xi \frac{ |\xi|^{2-2\rho}}{\la \xi \ra^{2\rho}} \, d\xi \right)^{1/2}
\end{equation}
Carrying out the $\xi$ integral over the region $|\xi|\leq |\tau_1-\xi_1^3|^{1/3}$ gives  
\begin{equation*} 
\int_\xi \frac{ |\xi|^{2-2\rho}}{\la \xi \ra^{2\rho}} \, d\xi \leq \la \tau_1-\xi_1^3 \ra^{1-\frac{4}{3}\rho}
\end{equation*}
and thus \eqref{AE:207} is bounded by
$$
 \la \tau_1-\xi_1^3 \ra^{1+\frac{1}{3}\rho-3b}
$$
which is bounded provided $b\geq \frac{5}{12}$.\\
\textit{Case 2}.  $3|\xi\xi_1\xi_2| \geq \frac{1}{2} |\tau_1-\xi_1^3|$. \\
In this case, $|\xi|\leq \frac{1}{10}|\xi_1|$.  Indeed, if $|\xi_1|\leq 10|\xi|$, then $3|\xi\xi_1\xi_2| \leq 330 |\xi|^3 \leq \frac{1}{3}|\tau_1-\xi_1^3|$.  Let $u=\tau_1-\xi_1^3-3\xi_1(\xi-\xi_1)\xi + \xi^3$, $du=3\xi_1(-2\xi+\xi_1)+3\xi^2$.  Now $3|\xi|^2 \leq \frac{3}{100}|\xi_1|^2$ and $3|\xi_1(-2\xi+\xi_1)| \geq \frac{12}{5}|\xi_1|^2$, and thus $3|\xi|^2 << 3|\xi_1(2\xi-\xi_1)|$.  We see that \eqref{AE:210} is bounded by
$$
 \frac{|\xi_1|^\rho}{\la \tau_1 -\xi_1^3 \ra^b} \left( \int_\xi \frac{|\xi|^2 |\xi_2|^{2\rho} |3 \xi_1(\xi_1-2\xi) + 3\xi^2|}{\la \xi \ra^{2\rho} \la \tau_1-\xi_1^3 - 3\xi\xi_1\xi_2 + \xi^3 \ra^{4b-1} |\xi_1(\xi_1-2\xi)|} \, d\xi \right)^{1/2}
$$
Using $|\xi_2|\sim |\xi_1|$, and $|\xi_1(\xi_1-2\xi)| \sim |\xi_1|^2$, this is controlled by
$$ \frac{|\xi_1|^{2\rho-1}}{\la \tau_1-\xi_1^3 \ra^b} \left( \int_{|u|\leq |\tau_1-\xi_1^3|} \frac{|\xi|^{2-2\rho}}{\la u \ra^{4b-1}} \, du \right)^{1/2}
$$
Using that $|\xi| \leq \dfrac{|\tau_1-\xi_1^3|}{|\xi_1|^2}$, this is controlled by
$$ \frac{|\xi_1|^{2\rho-1}|\tau_1-\xi_1^3|^{1-\rho}}{|\xi_1|^{2(1-\rho)}\la \tau_1-\xi_1^3 \ra^b} \left( \int_{|u|\leq |\tau_1-\xi_1^3|} \frac{du}{\la u \ra^{4b-1}} \right)^{1/2}
$$
Carrying out the $u$ integral, this is bounded by
$$ \frac{|\xi_1|^{4\rho-3}}{\la \tau_1-\xi_1^3 \ra^{\rho+3b-2}}
$$
which is bounded provided $b\geq \frac{2}{3}-\frac{1}{3}\rho$.\\

Now we address the range $\frac{3}{2}<s<3$.
It suffices to show
$$\iint_\ast \frac{|\xi| \la \tau \ra^{s/3} \hat{d}(\xi,\tau) }{\la \tau -\xi^3 \ra^b} \frac{\hat{g}_1(\xi_1,\tau_1)}{\la \tau_1-\xi_1^3 \ra^b \la \xi_1 \ra^s} \frac{\hat{g}_2(\xi_2,\tau_2)}{\la \tau_2 -\xi_2^3 \ra^b \la \xi_2 \ra^s}  \leq c\|d\|_{L^2} \|g_1\|_{L^2} \|g_2\|_{L^2}$$
for $\hat{d}\geq 0$, $\hat{g_1}\geq 0$, $\hat{g}_2\geq 0$, where $\ast$ indicates integration over $\xi$, $\xi_1$, $\xi_2$, subject to the constraint $\xi=\xi_1+\xi_2$, and over $\tau$, $\tau_1$, $\tau_2$, subject to the constraint $\tau=\tau_1+\tau_2$, under the assumption $|\tau|>>|\xi|^3$, since, for $s>0$ in the region $|\tau|\leq 2|\xi|^3$, $\|\partial_x(uv)\|_{Y_{s,-b}} \leq c\|\partial_x(uv)\|_{X_{s,-b}}$.  We shall show
\begin{equation} \label{AE:2200}
\frac{|\xi| \la \tau \ra^{s/3}}{\la \tau-\xi^3 \ra^b} \left( \int_{\xi_1} \int_{\tau_1} \frac{ d\xi_1 \, d\tau_1}{\la \tau_1-\xi_1^3 \ra^{2b} \la \xi_1 \ra^{2s} \la \tau_2-\xi_2^3 \ra^{2b} \la \xi_2 \ra^{2s}} \right)^{1/2} \leq c
\end{equation}
Since $\tau-\xi^3 +3\xi\xi_1\xi_2 = (\tau_2-\xi_2^3)+(\tau_1-\xi_1^3)$, by Lemma \ref{calc1}, we have that \eqref{AE:2200} is bounded by
\begin{equation} \label{AE:2201}
|\xi|\la \tau \ra^{\frac{s}{3}-b} \left( \int_{\xi_1} \frac{1}{\la \xi_1 \ra^{2s} \la \xi_2 \ra^{2s}} \frac{1}{\la \tau-\xi^3 + 3 \xi\xi_1\xi_2 \ra^{4b-1}} \, d\xi_1 \right)^{1/2}
\end{equation}
\noindent \textit{Case 1}. $|\xi_1|<< |\xi_2|$ or $|\xi_2|<<|\xi_1|$.  In this case, $3|\xi\xi_1\xi_2| << |\xi|^3$, which combined with $|\xi|^3<<|\tau|$, implies $\la \tau-\xi^3 + 3\xi\xi_1\xi_2 \ra \sim \la \tau \ra$.  Thus
\begin{align} 
\eqref{AE:2201} &\leq \la \tau \ra^{\frac{s}{3}-3b+\frac{1}{2}} \left( \int_{\xi_1} \frac{|\xi|^2 \, d\xi_1}{\la \xi_1 \ra^{2s} \la \xi_2 \ra^{2s}} \right)^{1/2} \notag \\
&\leq \la \tau \ra^{\frac{s}{3}-3b+\frac{1}{2}} \left( \int_{\xi_1} \frac{ d\xi_1}{\la \xi_1 \ra^{2s-2} \la \xi_2 \ra^{2s-2}} \right)^{1/2} \label{AE:2202}
\end{align}
Provided $s>\frac{3}{2}$ and $b>\frac{1}{9}s+\frac{1}{6}$, \eqref{AE:2202} is bounded.\\
\noindent \textit{Case 2}. $|\xi_1|\sim |\xi_2|$.  \\
\noindent \textit{Case 2A}.  $3|\xi\xi_1\xi_2|\sim |\tau|$ or $3|\xi\xi_1\xi_2| >> |\tau|$.  Then we ignore $\la \tau -\xi^3 +3\xi\xi_1\xi_2 \ra^{4b-1}$ in \eqref{AE:2201} and bound as:
$$\left( \int_{\xi_1} \frac{ |\xi|^2 \la \tau \ra^{\frac{2s}{3}-2b}}{\la \xi_1 \ra^{2s} \la \xi_2 \ra^{2s}} \, d\xi_1 \right)^{1/2}
$$
Using that $\la \tau \ra \leq c \la \xi \ra \la \xi_1 \ra \la \xi_2 \ra$, $\la \xi \ra \leq \la \xi_1 \ra+\la \xi_2 \ra$, and $\la \xi_1 \ra \sim \la \xi_2 \ra$, this is controlled by
$$\left( \int \frac{1}{\la \xi_1 \ra^{2s+6b-2}} \, d\xi_1 \right)^{1/2}
$$
Thus, we need $2s+6b-2>1$, which is automatically satisfied if $s>\frac{3}{2}$ and $b>0$.\\
\noindent \textit{Case 2B}.  $3|\xi\xi_1\xi_2| << |\tau|$.  Here, we just follow the method of Case 1. 

Thus we have estimate \eqref{AE:3151} for $-\frac{3}{4}<s<-\frac{1}{2}$, and $\frac{3}{2}<s<3$.  The result in the full range $-\frac{3}{4}<s<3$ follows by interpolation.
\end{proof}

\section{The left half-line problem}
\label{S:left}

We now carry out the proof of Theorem \ref{T:main}\eqref{I:left}.  We first return to the linearized version of \eqref{SE:101}.  Consider $-1< \lambda_1, \lambda_2 < 1$, $h_1,h_2 \in C_0^\infty(\mathbb{R}^+)$, and let
$$u(x,t) = \mathcal{L}_-^{\lambda_1}h_1(x,t) + \mathcal{L}_-^{\lambda_2}h_2(x,t)$$
By Lemma \ref{L:decayest}, $u(x,t)$ is continuous in $x$ at $x=0$ and by Lemma \ref{L:valueatzero},
$$u(0,t) = 2\sin(\tfrac{\pi}{3}\lambda_1 + \tfrac{\pi}{6})h_1(t) + 2 \sin(\tfrac{\pi}{3}\lambda_2+\tfrac{\pi}{6})h_2(t)$$
By the definition \eqref{AE:2020},
$$\partial_x u(x,t) = \mathcal{L}_-^{\lambda_1-1}\mathcal{I}_{-1/3}h_1(x,t) + \mathcal{L}_-^{\lambda_2-1}\mathcal{I}_{-1/3}h_2(x,t)$$
By Lemma \ref{L:decayest}, $\partial_xu(x,t)$ is continuous in $x$ at $x=0$ and by Lemma \ref{L:valueatzero},
$$\partial_xu(0,t) = 2\sin(\tfrac{\pi}{3}\lambda_1 - \tfrac{\pi}{6})h_1(t) + 2 \sin(\tfrac{\pi}{3}\lambda_2-\tfrac{\pi}{6})h_2(t)$$
Combining,
$$\begin{bmatrix} u(0,t) \\ \mathcal{I}_{-1/3}[\partial_x u(0,\cdot)](t) \end{bmatrix} = 2 \begin{bmatrix} \sin(\tfrac{\pi}{3}\lambda_1 + \tfrac{\pi}{6}) & \sin(\tfrac{\pi}{3}\lambda_2 + \tfrac{\pi}{6}) \\ \sin(\tfrac{\pi}{3}\lambda_1 - \tfrac{\pi}{6}) & \sin(\tfrac{\pi}{3}\lambda_2 - \tfrac{\pi}{6}) \end{bmatrix} \begin{bmatrix} h_1(t) \\ h_2(t) \end{bmatrix}$$
By basic trigonometric identities, this $2\times 2$ matrix has determinant $\sqrt{3}\sin \frac{\pi}{3}(\lambda_2-\lambda_1)$ which is $\neq 0$ provided $\lambda_1-\lambda_2 \neq 3n$ for $n\in \mathbb{Z}$.  Thus, for any $-1<\lambda_1, \lambda_2 <1$, with $\lambda_1\neq \lambda_2$, if we are given $g_1(t)$, $g_2(t)$ and we set
$$
\begin{bmatrix}
h_1(t) \\ h_2(t)
\end{bmatrix}
= A
\begin{bmatrix}
g_1(t) \\ \mathcal{I}_{1/3}g_2(t) 
\end{bmatrix}
$$
where
$$A= \frac{1}{2\sqrt 3\sin [ \frac{\pi}{3}(\lambda_2-\lambda_1)]} 
\begin{bmatrix} 
\sin(\tfrac{\pi}{3}\lambda_2 - \tfrac{\pi}{6}) & -\sin(\tfrac{\pi}{3}\lambda_2 + \tfrac{\pi}{6}) \\ 
-\sin(\tfrac{\pi}{3}\lambda_1 - \tfrac{\pi}{6}) & \sin(\tfrac{\pi}{3}\lambda_1 + \tfrac{\pi}{6}) 
\end{bmatrix}$$
then $u(x,t)$ solves
$$
\left\{
\begin{aligned}
& \partial_t u + \partial_x^3u = 0 && \text{for }x<0 \\
&u(x,0) = 0 \\
&u(0,t) = g_1(t) \\
&\partial_x u(0,t)=g_2(t)
\end{aligned}
\right.
$$
If we take $-1<\lambda_1, \lambda_2 <1$, $\lambda_1 \neq \lambda_2$, and set
$$\Lambda w(x,t) = \theta(t) e^{-t\partial_x^3}\phi(x) - \tfrac{1}{2}\theta(t) \mathcal{D}\partial_x w^2(x,t) + \theta(t)\mathcal{L}_-^{\lambda_1}h_1(x,t) + \theta(t)\mathcal{L}_-^{\lambda_2}h_2(x,t)$$
where
$$
\begin{bmatrix}
h_1(t) \\ h_2(t) 
\end{bmatrix}
= A
\begin{bmatrix}
g_1(t) - \theta(t) e^{-t\partial_x^3}\phi|_{x=0} + \tfrac{1}{2}\theta(t)\mathcal{D}\partial_xw^2(0,t) \\
\theta(t)\mathcal{I}_{1/3}( g_2 - \theta \partial_xe^{-\cdot \partial_x^3}\phi|_{x=0} + \tfrac{1}{2}\theta \partial_x\mathcal{D}\partial_xw^2(0,\cdot))(t)
\end{bmatrix}
$$
Then $(\partial_t + \partial_x^3)\Lambda w(x,t) = -\frac{1}{2}\partial_x w^2(x,t)$ for $x<0$, $0<t<1$, in the sense of distributions.  
We have
\begin{equation}
\label{E:351}
\begin{aligned}
\indentalign \|h_1\|_{H_0^\frac{s+1}{3}(\mathbb{R}^+)} + \|h_2\|_{H_0^\frac{s+1}{3}(\mathbb{R}^+)} \\
&\leq 
\begin{aligned}[t]
&c\|g_1(t) -\theta(t)e^{-t\partial_x^3}\phi|_{x=0} + \tfrac{1}{2}\theta(t)\mathcal{D}\partial_xw^2(0,t)\|_{H_0^\frac{s+1}{3}(\mathbb{R}^+)} \\
&+c\|\theta(t)\mathcal{I}_{1/3}( g_2 - \theta \partial_xe^{-\cdot \partial_x^3}\phi|_{x=0} + \tfrac{1}{2}\theta \partial_x\mathcal{D}\partial_xw^2(0,\cdot))(t)\|_{H_0^\frac{s+1}{3}(\mathbb{R}^+)}
\end{aligned}
\end{aligned}
\end{equation}
By Lemma \ref{L:G}\eqref{I:Gtt}, $\|g_1(t)-\theta(t)e^{-t\partial_x^3}\phi|_{x=0}\|_{H_t^\frac{s+1}{3}}\leq c\|g_1\|_{H^\frac{s+1}{3}}+c\|\phi\|_{H^s}$.  If $-\frac{3}{4}<s<\frac{1}{2}$, then $\frac{1}{12}<\frac{s+1}{3}<\frac{1}{2}$, and Lemma \ref{JK35} shows that $g_1(t)-\theta(t)e^{-t\partial_x^3}\phi|_{x=0}\in H_0^\frac{s+1}{3}(\mathbb{R}_t^+)$ with comparable norm.  If $\frac{1}{2}<s<\frac{3}{2}$, then $\frac{1}{2}<\frac{s+1}{3}<\frac{5}{6}$ and by the compatibility condition, $g_1(t)-\theta(t)e^{-t\partial_x^3}\phi|_{x=0}$ has a well-defined value of $0$ at $t=0$.  By Lemma \ref{JK37}, $g_1(t)-\theta(t)e^{-t\partial_x^3}\phi|_{x=0}$ also belongs to $H_0^\frac{s+1}{3}(\mathbb{R}_t^+)$ with comparable norm.  The conclusion then, is that if $-\frac{3}{4}<s<\frac{3}{2}$, $s\neq \frac{1}{2}$, then
$$\|g_1(t)-\theta(t)e^{-t\partial_x^3}\phi|_{x=0}\|_{H_0^\frac{s+1}{3}(\mathbb{R}^+)} \leq c\|g_1\|_{H^\frac{s+1}{3}} + c\|\phi\|_{H^s}$$
By Lemmas \ref{L:Di}\eqref{I:Ditt}, \ref{L:bilinear},
$$\| \theta(t)\mathcal{D}\partial_xw^2(0,t) \|_{H_0^\frac{s+1}{3}(\mathbb{R}_t^+)} \leq c\|w\|_{X_{s,b}\cap D_\alpha}^2$$
By Lemma \ref{L:G}\eqref{I:Gdtt}, $\|g_2(t) - \theta(t)\partial_xe^{-t\partial_x^3}\phi|_{x=0}\|_{H_t^{s/3}} \leq c\|g_2\|_{H^{s/3}} + c\|\phi\|_{H^s}$.  If $-\frac{3}{4}<s<\frac{3}{2}$, then $\frac{s}{3}<\frac{1}{2}$ and by Lemma \ref{JK35}, $g_2(t) - \theta(t)\partial_xe^{-t\partial_x^3}\phi|_{x=0}\in H_0^{s/3}(\mathbb{R}^+)$ with comparable norm.  By Lemma \ref{L:RL2}, 
$$\| \theta(t)\mathcal{I}_{1/3}( g_2 - \theta \partial_xe^{-\cdot \partial_x^3}\phi|_{x=0}) \|_{H_0^\frac{s+1}{3}(\mathbb{R}^+)} \leq c\|g_1\|_{H^\frac{s+1}{3}} + c\|\phi\|_{H^s}$$
By Lemmas \ref{L:RL2}, \ref{L:Di}\eqref{I:Didtt}, \ref{L:bilinear},
$$\| \theta(t)\mathcal{I}_{1/3}(\theta \partial_x\mathcal{D}\partial_xw^2(0,\cdot))(t) \|_{H_0^\frac{s+1}{3}(\mathbb{R}_t^+)} \leq c\|w\|_{X_{s,b}\cap D_\alpha}^2$$
Combining the above estimates with \eqref{E:351}, we obtain
\begin{equation}
\label{E:350}
\|h_1\|_{H_0^{\frac{s+1}{3}}(\mathbb{R}^+)} + \|h_2\|_{H_0^{\frac{s+1}{3}}(\mathbb{R}^+)} \leq c\|g_1\|_{H_t^\frac{s+1}{3}}+ \|g_2\|_{H_t^{s/3}} + c\|\phi\|_{H^s} + c\|w\|_{X_{s,b}\cap D_\alpha}^2
\end{equation}
By Lemmas \ref{L:G}\eqref{I:Gst}, \ref{L:Di}\eqref{I:Dist}, \ref{L:Dbf}\eqref{I:Dbfst}, \ref{L:bilinear}, and \eqref{E:350}
$$\| \Lambda w(x,t) \|_{C(\mathbb{R}_t; H_x^s)} \leq c\|\phi\|_{H^s} + c\|g_1\|_{H^\frac{s+1}{3}} + c\|g_2\|_{H^\frac{s}{3}}+ c\|w\|_{X_{s,b}\cap D_\alpha}^2$$
provided $b(s)\leq b< \frac{1}{2}$ (where $b(s)$ is specified by Lemma \ref{L:bilinear}), $s-\frac{5}{2}<\lambda_1<s+\frac{1}{2}$, $s-\frac{5}{2}<\lambda_2<s+\frac{1}{2}$, $\alpha>\frac{1}{2}$.  In the sense of $C(\mathbb{R}_t; H_x^s)$, $w(x,0)=\phi(x)$.  By Lemmas \ref{L:G} \eqref{I:Gtt}, \ref{L:Di}\eqref{I:Ditt}, \ref{L:Dbf}\eqref{I:Dbftt}, \ref{L:bilinear}, and \eqref{E:350}
$$\|\Lambda w(x,t) \|_{C(\mathbb{R}_x;H_t^\frac{s+1}{3})} \leq c\|\phi\|_{H^s} + c\|g_1\|_{H^\frac{s+1}{3}} + c\|g_2\|_{H^\frac{s}{3}}+ c\|w\|_{X_{s,b}\cap D_\alpha}^2$$
provided $b(s)<b<\frac{1}{2}$.  In the sense of $C(\mathbb{R}_x; H_t^\frac{s+1}{3})$, $\Lambda w(0,t) = g_1(t)$ for $0\leq t\leq 1$.  By Lemmas \ref{L:G}\eqref{I:Gdtt}, \ref{L:Di}\eqref{I:Didtt}, \ref{L:Dbf}\eqref{I:Dbfdtt}, \ref{L:bilinear}, and \eqref{E:350}
$$\|\partial_x \Lambda w(x,t) \|_{C(\mathbb{R}_x;H_t^\frac{s}{3})} \leq c\|\phi\|_{H^s} + c\|g_1\|_{H^\frac{s+1}{3}} + c\|g_2\|_{H^\frac{s}{3}}+ c\|w\|_{X_{s,b}\cap D_\alpha}^2$$
provided $b(s)<b<\frac{1}{2}$, and in the sense of $C(\mathbb{R}_x;H^{s/3}_t)$, $\partial_xw(0,t)=g_2(t)$ for $0\leq t\leq 1$.  By Lemma \ref{L:G}\eqref{I:GBse}, \ref{L:Di}\eqref{I:DiBse}, \ref{L:Dbf}\eqref{I:DbfBse}, \ref{L:bilinear}, and \eqref{E:350}, we have
$$\|\Lambda w \|_{X_{s,b}\cap D_\alpha} \leq c\|\phi\|_{H^s} + c\|g_1\|_{H^\frac{s+1}{3}} + c\|g_2\|_{H^\frac{s}{3}}+ c\|w\|_{X_{s,b}\cap D_\alpha}^2$$
provided $s-1\leq \lambda_1<s+\frac{1}{2}$, $s-1\leq \lambda_2 < s+\frac{1}{2}$, $\lambda_1<\frac{1}{2}$, $\lambda_2<\frac{1}{2}$, $\alpha\leq \frac{s-\lambda_1+2}{3}$, $\alpha\leq \frac{s-\lambda_2+2}{3}$, $b(s)<b<\frac{1}{2}$, and $\frac{1}{2}<\alpha\leq 1-b$.
  
Collectively, the restrictions are $-\frac{3}{4}<s<\frac{3}{2}$, $s\neq \frac{1}{2}$, $b(s)<b<\frac{1}{2}$, 
\begin{equation} 
\label{E:360}
\begin{aligned}
&s-1\leq \lambda_1 <s+\tfrac{1}{2} &\qquad & -1< \lambda_1<\tfrac{1}{2}\\
& s-1\leq \lambda_2 < s+\tfrac{1}{2} && -1<\lambda_2 < \tfrac{1}{2}
\end{aligned}
\end{equation}
\begin{equation}
\label{E:361}
\begin{aligned}
&\tfrac{1}{2}< \alpha \leq \tfrac{s-\lambda_1+2}{3} \\
&\tfrac{1}{2}< \alpha \leq \tfrac{s-\lambda_2+2}{3} \\
&\alpha \leq 1-b
\end{aligned}
\end{equation}
Since $s<\frac{3}{2} \Longrightarrow s-1<\frac{1}{2}$ and $s>-\frac{3}{4} \Longrightarrow s+\frac{1}{2}> -\frac{1}{4}$, and thus we can find $\lambda_1\neq \lambda_2$ meeting the restriction \eqref{E:360}.  (Note that for $s<-\frac{1}{2}$, we cannot use $\lambda=0$, the operator used in \cite{CK02}).  The conditions $\lambda_1<s+\frac{1}{2}$, $\lambda_2<s+\frac{1}{2}$ imply that $\frac{s-\lambda_1+2}{3}>\frac{1}{2}$, $\frac{s-\lambda_2+2}{3}>\frac{1}{2}$, and thus we can meet the requirements expressed in \eqref{E:360}.

Define a space $Z$ by the norm
$$\|w\|_Z = \|w\|_{C(\mathbb{R}_t;H_x^s)} + \|w\|_{C(\mathbb{R}_x;H_t^\frac{s+1}{3})} +  \|\partial_x w\|_{C(\mathbb{R}_x;H_t^\frac{s+1}{3})} + \|w\|_{X_{s,b}\cap D_\alpha}$$
By the above estimates
$$\|\Lambda w \|_Z \leq   c\|\phi\|_{H^s} + c\|g_1\|_{H^\frac{s+1}{3}} + c\|g_2\|_{H^\frac{s}{3}}+ c\|w\|_Z^2$$
Now
\begin{align*}
\indentalign \Lambda w_1(x,t) - \Lambda w_2(x,t) \\
&= 
\begin{aligned}[t]
&-\tfrac{1}{2}\theta(t)\mathcal{D}\partial_x(w_1-w_2)(w_1+w_2)(x,t) + \theta(t)\mathcal{L}_-^{\lambda_1}h_1(x,t) \\
&+ \theta(t)\mathcal{L}_-^{\lambda_2}h_2(x,t)
\end{aligned}
\end{align*}
where
$$
\begin{bmatrix}
h_1(t) \\ h_2(t)
\end{bmatrix} 
=\tfrac{1}{2}A
\begin{bmatrix}
\theta(t) \mathcal{D}\partial_x(w_1-w_2)(w_1+w_2)(0,t) \\
\theta(t) \mathcal{I}_{1/3}( \theta \partial_x \mathcal{D} \partial_x(w_1-w_2)(w_1+w_2)(0,\cdot))(t)
\end{bmatrix}
$$
By similar arguments, we can show
$$\|\Lambda w_1- \Lambda w_2 \|_2 \leq c(\|w_1\|_Z+\|w_2\|_Z)(\|w_1-w_2\|_Z)$$
By taking $\|\phi\|_{H^s} + \|g_1\|_{H^\frac{s+1}{3}} + \|g_2\|_{H^\frac{s}{3}}\leq \delta$ for $\delta>0$ suitably small, we obtain a fixed point ($\Lambda u =u$) in $Z$.

Theorem \ref{T:main}\eqref{I:left} follows by the standard scaling argument. Suppose we are given data $\tilde \phi$, $\tilde g_1$, and $\tilde g_2$ of arbitrary size for the problem \eqref{SE:101}, and we seek a solution $\tilde u$.  For $0 \leq \lambda \ll 1$ (to be selected in a moment) set $\phi(x) = \lambda^2\tilde \phi(x)$, $g_1(t) = \lambda^2\tilde g_1(t)$, $g_2(t) = \lambda^3 \tilde g_2(\lambda^3 t)$.  Take $\lambda$ sufficiently small so that 
\begin{align*}
\indentalign \|\phi\|_{H^s} + \|g_1\|_{H^\frac{s+1}{3}}+ \|g_2\|_{H^\frac{s}{3}}  \\
&\leq \lambda^\frac{3}{2}\la \lambda^s \ra \|\tilde \phi \|_{H^s} + \lambda^\frac{1}{2} \la \lambda \ra^{s+1}\| \tilde g_1 \|_{H^\frac{s+1}{3}} + \lambda^\frac{3}{2} \la \lambda \ra^{s}\| \tilde g_2 \|_{H^\frac{s}{3}}\\
&\leq \delta
\end{align*}
By the above argument, there is a solution $u(x,t)$ on $0\leq t\leq 1$.  Then $\tilde u(x,t) = \lambda^{-2} u( \lambda^{-1} x, \lambda^{-3} t)$ is the desired solution on $0\leq t \leq \lambda^3$.

\section{The right half-line problem}
\label{S:right}

Now we prove Theorem \ref{T:main}\eqref{I:right}.  Suppose $-1<\lambda<1$ and we are given $f\in C_0^\infty(\mathbb{R}^+)$.  Let $u(x,t)=e^{-\pi \lambda i}\mathcal{L}_+^\lambda f(x,t)$.  Then by Lemma \ref{L:decayest}, $u(x,t)$ is continuous in $x$ at $x=0$ and by Lemma \ref{L:valueatzero}, $u(0,t)=f(t)$.  Then $u(x,t)$ solves
$$
\left\{
\begin{aligned}
&\partial_t u +\partial_x^3 u = 0 \\
&u(x,0) = 0 \\
&u(0,t) = f(t) 
\end{aligned}
\right.
$$
Therefore, to address the nonlinear problem \eqref{SE:101} with given data $f$ and $\phi$, take $-1<\lambda<1$ and set
$$\Lambda w(x,t) = \theta(t) e^{-t\partial_x^3}\phi(x) - \tfrac{1}{2}\theta(t)\mathcal{D}\partial_x w^2(x,t) + \theta(t)\mathcal{L}_+^\lambda h(x,t)$$
where
$$h(t) = e^{-\pi i\lambda}[ f(t) - \theta(t)e^{-t\partial_x^3}\phi|_{x=0} + \tfrac{1}{2}\theta(t)\mathcal{D}\partial_xw^2(0,t)]$$
Then 
$$(\partial_t + \partial_x^3)\Lambda w(x,t) = -\frac{1}{2}\partial_x w^2(x,t).$$  
By Lemma \ref{L:G}\eqref{I:Gtt}, $\|f(t)-\theta(t)e^{-t\partial_x^3}\phi|_{x=0} \|_{H^\frac{s+1}{3}} \leq c\|f\|_{H^\frac{s+1}{3}}+c\|\phi\|_{H^s}$.  If $-\frac{3}{4}<s<\frac{1}{2}$, then $\frac{1}{12}<\frac{s+1}{3}<\frac{1}{2}$ and Lemma \ref{JK35} shows that $f(t)-\theta(t)e^{-t\partial_x^3}\phi|_{x=0}\in H_0^\frac{s+1}{3}$ with comparable norm.  If $\frac{1}{2}<s<\frac{3}{2}$, then $\frac{1}{2}<\frac{s+1}{3}<\frac{5}{6}$ and by the compatibility condition, $f(t)-\theta(t)e^{-t\partial_x^3}\phi|_{x=0}$ has a well-defined value of $0$ at $t=0$.  By Lemma \ref{JK37}, $f(t)-\theta(t)e^{-t\partial_x^3}\phi|_{x=0}\in H_0^\frac{s+1}{3}(\mathbb{R}^+)$ with comparable norm.  The conclusion, then, is that if $-\frac{3}{4}<s<\frac{3}{2}$, $s\neq \frac{1}{2}$, then
$$\|f(t)-\theta(t)e^{-t\partial_x^3}\phi|_{x=0}\|_{H^\frac{s+1}{3}(\mathbb{R}^+)} \leq c\|f\|_{H^\frac{s+1}{3}} + c\|\phi\|_{H^s}$$
By Lemma \ref{L:Di}\eqref{I:Ditt}, \ref{L:bilinear},
$$\|\theta(t)\mathcal{D}\partial_x w^2(0,t)\|_{H_0^\frac{s+1}{3}(\mathbb{R}^+)} \leq c\|w\|_{X_{s,b}\cap D_\alpha}^2$$
Combining, we obtain
\begin{equation}
\label{E:362}
\|h\|_{H_0^\frac{s+1}{3}(\mathbb{R}^+)} \leq c\|f\|_{H^\frac{s+1}{3}} + c\|\phi\|_{H^s} + c\|w\|_{X_{s,b}\cap D_\alpha}^2
\end{equation}

We then proceed in the manner of \S \ref{S:left} to complete the proof of Theorem \ref{T:main}\eqref{I:right}.

\section{The line segment problem}
\label{S:linesegment}

We now turn to the line segment problem \eqref{SE:102}.  By the standard scaling argument, it suffices to show that $\exists \; \delta>0$ and $\exists \; L_1>>0$ such that for any $L>L_1$ and data $f$, $g_1$, $g_2$, $\phi$ satisfying
$$\|f\|_{H^\frac{s+1}{3}(\mathbb{R}^+)} + \|g_1\|_{H^\frac{s+1}{3}(\mathbb{R}^+)} + \|g_2\|_{H^\frac{s}{3}(\mathbb{R}^+)}+\|\phi\|_{H^s(0,L)} \leq \delta$$
we can solve \eqref{SE:102} with $T=1$.  By the techniques employed in the previous two sections, it suffices to show that for all boundary data $f$, $g_1$, $g_2$, there exists $u$ solving the linear problem
\begin{equation}
\label{E:366}
\left\{
\begin{aligned}
&\partial_tu + \partial_x^3 u = 0  && \text{for }(x,t)\in (0,L)\times (0,1) \\
& u(0,t) = f(t) && \text{for }t\in (0,1)\\
& u(L,t) = g_1(t) && \text{for }t\in (0,1)\\
&\partial_xu(L,t) = g_2(t) && \text{for }t\in (0,1)\\
&u(x,0) = 0 && \text{for }x\in (0,L)
\end{aligned}
\right.
\end{equation}
such that
\begin{equation}
\label{E:367}
\begin{aligned}
\indentalign \|u\|_{C(\mathbb{R}_t; H^s_x)} + \| u\|_{C(\mathbb{R}_x; H^\frac{s+1}{3}_t)} + \|\partial_x u \|_{C(\mathbb{R}_x; H^\frac{s}{3}_t)}+ \|u\|_{X_{s,b}\cap D_\alpha} \\
&\leq  \|f\|_{H^\frac{s+1}{3}(\mathbb{R}^+)} + \|g_1\|_{H^\frac{s+1}{3}(\mathbb{R}^+)} + \|g_2\|_{H^\frac{s}{3}(\mathbb{R}^+)}
\end{aligned}
\end{equation}
Let
\begin{align*}
\mathcal{L}_1h_1(x,t)&= \mathcal{L}_-^{\lambda_1}h_1(x-L,t) \\
\mathcal{L}_2h_2(x,t)&= \mathcal{L}_-^{\lambda_2}h_2(x-L,t) \\
\mathcal{L}_3h_3(x,t)&=\mathcal{L}_+^{\lambda_3}h_3(x,t)
\end{align*}
By Lemma \ref{L:valueatzero} and the estimates in \S \ref{S:estimates}, solving \eqref{E:366}, \eqref{E:367} amounts to showing that the matrix equation
\begin{equation}
\label{E:364}
(g_1,\mathcal{I}_{1/3}g_2, f)^T = (E_L+K_L)(h_1,h_2,h_3)^T
\end{equation}
has a bounded inverse, where 
$$E_L=\begin{bmatrix}
2\sin( \frac{\pi}{3}\lambda_1+\frac{\pi}{6})  &&
2\sin( \frac{\pi}{3}\lambda_2+\frac{\pi}{6})  &&
 0 \\
2\sin( \frac{\pi}{3}\lambda_1-\frac{\pi}{6})  &&
2\sin( \frac{\pi}{3}\lambda_2-\frac{\pi}{6})  &&
 0 \\
\mathcal{L}_1\big|_{x=0} && \mathcal{L}_2\big|_{x=0} && e^{i\pi \lambda_3}
\end{bmatrix},$$
$$
K_L=\begin{bmatrix}
0 && 0 && \mathcal{L}_3\big|_{x=L} \\
0 && 0 && \mathcal{I}^t_{1/3}(\partial_x\mathcal{L}_3)\big|_{x=L} \\
0 && 0 && 0
\end{bmatrix}$$
The matrix operator $E_L$ is invertible with inverse 
$$E_L^{-1}= \begin{bmatrix}
\dfrac{\sin(\frac{\pi}{3}\lambda_2-\frac{\pi}{6})}{\sqrt{3}\sin(\frac{\pi}{3}\lambda_2-\frac{\pi}{3}\lambda_1)} &&
\dfrac{-\sin(\frac{\pi}{3}\lambda_2+\frac{\pi}{6})}{\sqrt{3}\sin(\frac{\pi}{3}\lambda_2-\frac{\pi}{3}\lambda_1)} &&
0 \\
\dfrac{-\sin(\frac{\pi}{3}\lambda_1-\frac{\pi}{6})}{\sqrt{3}\sin(\frac{\pi}{3}\lambda_2-\frac{\pi}{3}\lambda_1)} &&
\dfrac{\sin(\frac{\pi}{3}\lambda_1+\frac{\pi}{6})}{\sqrt{3}\sin(\frac{\pi}{3}\lambda_2-\frac{\pi}{3}\lambda_1)} &&
0 \\
A_1 && A_2 && e^{-i\pi\lambda_3}
\end{bmatrix}$$
where
$$A_1 =\frac{ \sqrt{3} e^{-i\pi\lambda_3}\sin( \frac{\pi}{3}\lambda_1-\frac{\pi}{6})}{\sin(\frac{\pi}{3}\lambda_2-\frac{\pi}{3}\lambda_1)} \mathcal{L}_2\big|_{x=0} - \frac{ \sqrt{3} e^{-i\pi\lambda_3} \sin( \frac{\pi}{3}\lambda_2-\frac{\pi}{6})}{\sin(\frac{\pi}{3}\lambda_2-\frac{\pi}{3}\lambda_1)}\mathcal{L}_1\big|_{x=0}$$
and
$$A_2 = \frac{ -\sqrt{3}e^{-i\pi\lambda_3} \sin( \frac{\pi}{3}\lambda_1+\frac{\pi}{6})}{\sin(\frac{\pi}{3}\lambda_2-\frac{\pi}{3}\lambda_1)} \mathcal{L}_2\big|_{x=0} + \frac{ \sqrt{3}e^{-i\pi\lambda_3} \sin( \frac{\pi}{3}\lambda_2+\frac{\pi}{6})}{\sin(\frac{\pi}{3}\lambda_2-\frac{\pi}{3}\lambda_1)}\mathcal{L}_1\big|_{x=0}$$
Since $\mathcal{L}_1\big|_{x=0}: H_0^\frac{s+1}{3}(\mathbb{R}^+)\to H_0^\frac{s+1}{3}(\mathbb{R}^+)$, $\mathcal{L}_2\big|_{x=0}: H_0^\frac{s+1}{3}(\mathbb{R}^+)\to H_0^\frac{s+1}{3}(\mathbb{R}^+)$ are bounded uniformly as $L\to +\infty$, the norm of $E_L^{-1}$ is uniformly bounded as $L\to +\infty$.  \eqref{E:364} becomes
\begin{equation} \label{E:363}
E_L^{-1}(g_1,\mathcal{I}_{1/3}g_2, f)^T = (I+E_L^{-1}K_L)(h_1,h_2,h_3)^T
\end{equation}
and we see that it suffices to show that $(I+E_L^{-1}K_L)$ is invertible.  We claim that $K_L:[ H_0^\frac{s+1}{3}(\mathbb{R}^+)]^3 \to [H_0^\frac{s+1}{3}(\mathbb{R}^+)]^3$ is bounded with norm $\to 0$ as $L\to +\infty$.  To show this, we need a refinement of Lemma \ref{L:Dbf}\eqref{I:Dbftt}. 

\begin{lemma}
For $-2<\lambda<1$ and $x>0$
$$\|\theta(t)\mathcal{L}_+^\lambda h(x,t) \|_{H_0^\frac{s+1}{3}(\mathbb{R}^+)} \leq c(x)\|h\|_{H_0^\frac{s+1}{3}(\mathbb{R}^+)}$$
where $c(x)\to 0$ as $x\to +\infty$.
\end{lemma}
\begin{proof}
$\mathcal{L}_+^\lambda f(x,t) = \mathcal{L}^0h(x,t)$ for $x>0$ by a uniqueness calculation.  By \eqref{E:Dbf}, 
\begin{align*}
\theta(t)\mathcal{L}^0h(x,t) &= \theta(t) \int_0^t \frac{\theta(2(t-t'))}{(t-t')^{1/3}} A\left( \frac{x}{(t-t')^{1/3}} \right) \mathcal{I}_{-2/3}h(t') \, dt' \\
&= -\theta(t) \int_0^t \partial_{t'}\left[ \frac{\theta(2(t-t'))}{(t-t')^{1/3}}  A\left( \frac{x}{(t-t')^{1/3}} \right)\right] \theta(4t')\mathcal{I}_{1/3}h(t') \, dt'
\end{align*}
Since $A(x)$ decay rapidly as $x\to +\infty$, we have
\begin{align*}
H(t) &:= -\partial_t \left[ \frac{\theta(2t)}{t^{1/3}} A \left( \frac{x}{t^{1/3}} \right) \chi_{t\geq 0} \right] \\
&= 
\begin{aligned}[t]
&-2\theta'(2t) x^{-1} \left( \frac{x}{t^{1/3}} \right) A\left( \frac{x}{t^{1/3}} \right) \chi_{t\geq 0} \\
&+ \tfrac{1}{3}\theta(2t) x^{-4} \left( \frac{x}{t^{1/3}} \right)^4 A\left( \frac{x}{t^{1/3}} \right) \chi_{t\geq 0} \\
&+ \tfrac{1}{3}\theta(2t) x^{-4}\left( \frac{x}{t^{1/3}} \right)^5 A' \left( \frac{x}{t^{1/3}} \right) \chi_{t\geq 0}
\end{aligned}
\end{align*}
so that $\mathcal{L}^0h(x,t) = \theta(t) H\ast (\theta(4\cdot)\mathcal{I}_{1/3}h)(t)$.  By the asymptotic properties of $A(x)$ as $x\to +\infty$,
$$\|\hat{H}\|_{L^\infty} \leq \|H\|_{L^1} \leq \sup_{x\geq \frac{x}{2}} ( |x^4A(x)| + |x^5A'(x)|)\to 0 \text{ as }x\to +\infty$$
and we have 
$$\| \mathcal{L}^0h(x,t)\|_{H^\frac{s+1}{3}} \leq \|\hat{H}\|_{L^\infty}\| \theta(4t)\mathcal{I}_{1/3}h(t)\|_{H^\frac{s+1}{3}} \leq c(x)\|h\|_{H^\frac{s+1}{3}}$$
with $c(x)\to 0$ as $x\to +\infty$.  
\end{proof}

From the lemma, $\mathcal{L}_3|_{x=L}: H_0^\frac{s+1}{3}(\mathbb{R}^+)\to H_0^\frac{s+1}{3}(\mathbb{R}^+)$ and $\mathcal{I}_{1/3}(\partial_x \mathcal{L}_3)|_{x=L}= \mathcal{I}_{1/3}(\mathcal{L}_+^{\lambda_3-1}\mathcal{I}_{-1/3})|_{x=L}: H_0^\frac{s+1}{3}(\mathbb{R}^+)\to H_0^\frac{s+1}{3}(\mathbb{R}^+)$ are bounded with norm $\to 0$ as $L\to +\infty$.  Thus $K_L:[ H_0^\frac{s+1}{3}(\mathbb{R}^+)]^3 \to [H_0^\frac{s+1}{3}(\mathbb{R}^+)]^3$ enjoys the same property and $(I+E_L^{-1}K_L)$ has bounded (uniformly in $L$ as $a\to +\infty$) inverse in \eqref{E:363}.

\def\cprime{$'$}
\providecommand{\bysame}{\leavevmode\hbox to3em{\hrulefill}\thinspace}
\providecommand{\MR}{\relax\ifhmode\unskip\space\fi MR }
\providecommand{\MRhref}[2]{%
  \href{http://www.ams.org/mathscinet-getitem?mr=#1}{#2}
}
\providecommand{\href}[2]{#2}

\end{document}